\documentclass[12pt,final,3p,times]{elsarticle}
\usepackage{amsthm,amsmath,amssymb,amsfonts,graphicx,graphics,latexsym,exscale,cmmib57,dsfont,bbm,amscd}
\usepackage{mathrsfs}

\usepackage{tikz-cd}
\swapnumbers
\newtheorem{sthm}{Theorem}[subsection]
\newtheorem{slem}[sthm]{Lemma}

\newtheorem{thm}{Theorem}[section]
\newtheorem{lem}[thm]{Lemma}
\newtheorem{cor}[thm]{Corollary}
\newtheorem*{1.3.1D}{\ref{1.3.1}D~Lemma}
\newtheorem*{2.2.2D}{\ref{2.2.2}D~Lemma}
\newtheorem*{2.2.2E}{\ref{2.2.2}E~Theorem}
\newtheorem*{xthm}{Theorem}
\theoremstyle{definition}

\newtheorem{rem}[thm]{Remark}
\newtheorem{srem}[sthm]{Remark}

\newtheorem{sq}[sthm]{Question}
\newtheorem{exa}[thm]{Example}

\newtheorem{se}[thm]{}
\newtheorem{sse}[sthm]{}
\newtheorem*{xnote}{Note}

\newtheorem*{2.2.2F}{\ref{2.2.2}F~Example}

\numberwithin{equation}{section}

\journal{Topol. Appl.}

\begin{document}

\begin{frontmatter}
\title{On semi-openness of fiber-onto extensions of minimal semiflows and quasi-separable maps}

\author{Xiongping Dai}
\ead{xpdai@nju.edu.cn}
\address{School of Mathematics, Nanjing University, Nanjing 210093, People's Republic of China}
\author{Li Feng}
%\ead{li.feng@asurams.edu}
\address{Department of Mathematics and Computer Science,
Albany State University,
Albany, Georgia 31705, USA}

\author{Congying Lv}
%\ead{lvcongying@smail.nju.edu.cn}
%\address{School of Mathematics, Nanjing University, Nanjing 210093, People's %Republic of China}
\author{Yuxuan Xie}
%\ead{201501005@smail.nju.edu.cn}
\address{School of Mathematics, Nanjing University, Nanjing 210093, People's Republic of China}

\author{In memory of Professor Zuoling Zhou (1938--2026)}
%%%%%%%%%%%%%%%%%%%%%%%%%%%%%%
\begin{abstract}
The purpose of this paper is to find conditions for a continuous onto map $\phi\colon X\rightarrow Y$ and its induced map $\phi_*\colon\mathcal{M}^1(X)\rightarrow\mathcal{M}^1(Y)$  to be semi-open, where $X$, $Y$ are compact Hausdorff spaces and $\mathcal{M}^1(X)$, $\mathcal{M}^1(Y)$ are their Borel probability spaces. For that, we mainly prove the following results by using the structure theory of extensions of semiflows and inverse limit techniques:
\begin{enumerate}[(1)]
\item If $\phi$ is an extension of minimal flows, then $\phi_*$ is semi-open.
\item If $\phi$ is a quasi-separable fiber-onto extension of minimal semiflows, then $\phi$ and $\phi_*$ are semi-open.
\item If $Y$ is metrizable, then $\phi$ is semi-open if and only if $\phi_*$ is semi-open.
\end{enumerate}
In addition, if $X,Y$ are left-topological groups, $X$ is Lindel\"{o}f quasi-regular, $Y$ is Baire and if $\phi$ is a locally closed continuous onto equivariant mapping, then $\phi$ is semi-open (This is a generalization of Pontryagin's open-mapping theorem). 
\end{abstract}

\end{frontmatter}
%%% ----------------------------------------------------------------------
%%% ----------------------------------------------------------------------

% Text of article.
\section{Introduction}
In this paper we mainly study the semi-openness of some canonical continuous onto mappings between topological spaces.
Let $\mathscr{O}(X)$, for any topological space $X$, stand for the family of open non-void subsets of $X$, $\mathfrak{N}_x(X)$ the filter of neighborhoods of $x$ in $X$; and let $\textrm{cl}\,A=\bar{A}$ and $\textrm{int}\,A$ be the closure and interior of $A$ in $X$ for any set $A\subset X$, respectively. If $X$ is a uniform/completely-regular space, then $\mathscr{U}_X$ denotes a compatible uniform structure for $X$ (cf.~\cite[Chap.~6]{K55}).

\subsection{Induced spaces and maps}
\begin{sse}[$f\colon X\rightarrow Y$ vs $2^f\colon 2^X\rightarrow2^Y$]
Let $2^X$ be the \textit{hyperspace} of closed non-void subsets of a topological space $X$
endowed with the {\it Vietoris topology}~(cf.~\cite[IX.1]{G76} or \cite[II.1]{Wo}), for which a base is given by the sets of the form
$$
\langle U_1,\dotsc,U_n\rangle=\left\{A\in 2^X\,|\,A\subseteq U_1\cup\dotsm\cup U_n\ \&\ A\cap U_i\not=\emptyset\ \forall i=1,\dotsc,n\right\},\ n\in\mathbb{N}\ \&\ U_i\in\mathscr{O}(X).
$$
Then, $X$ is compact Hausdorff if and only if so is $2^X$; $X$ is metrizable if and only if so is $2^X$ (cf.~\cite[Thm.~II.1.1]{Wo}). Clearly, $f\colon X\rightarrow Y$ is a continuous (onto) mapping between compact Hausdorff spaces if and only if so is its induced map $2^f\colon 2^X\rightarrow 2^Y$, $K\mapsto f[K]$ (cf.~\cite[Thm.~II.1.3a]{Wo}). 

It is well known that the topological properties of $2^X$ and $2^{2^X}$ are important for us to obtain information on the structure of the space $X$ (cf., e.g., \cite{K42, N68, IN}). In the case of $f\colon X\rightarrow X$ and $2^f\colon 2^X\rightarrow2^X$, some dynamics, such as ergodicity, recurrence, almost periodicity, equicontinuity, and disjointness, of $f$ and $2^f$ may be described via each other (cf., e.g., \cite{K75, G76, AG77, Wo, GW95, G00, CKN, LYY, AAN, LOYZ, HSY, N22, JO}).
\end{sse}

\begin{sse}[$f\colon X\rightarrow Y$ vs $f_*\colon \mathcal{M}^1(X)\rightarrow\mathcal{M}^1(Y)$]
Let $X$ be a compact Hausdorff space. By $\mathcal{M}^1(X)$ it means the set of all regular Borel probabilities on $X$, which is equipped with the weak* topology; that is, a net $\mu_\alpha\to\mu$ in $\mathcal{M}^1(X)$ iff $\mu_\alpha(\varphi)\to\mu(\varphi)$ for all $\varphi\in C(X)$, where $C(X)$ is the set of all continuous real-valued functions on $X$.
Then $\mathcal{M}^1(X)$ is a compact Hausdorff space; and moreover, $\mathcal{M}^1(X)\rightarrow2^X$, given by $\mu\mapsto\textrm{supp}(\mu)$, is a lower semi-continuous map (cf.~\cite[Lem.~VII.1.4]{Wo}). 

If $f\colon X\rightarrow Y$ is a continuous mapping between compact Hausdorff spaces, then its induced map $f_*\colon\mathcal{M}^1(X)\rightarrow\mathcal{M}^1(Y)$, ${\mu\mapsto \mu\circ f^{-1}}$ is continuous affine.
Let $\delta_x\in\mathcal{M}^1(X)$ be the Dirac measure at $x\in X$ and
$\delta[X]=\{\delta_x\,|\,x\in X\}$. If $f$ is onto, then by $f_*[\delta[X]]=\delta[Y]$ and $\mathcal{M}^1(Y)=\overline{\mathrm{co}}\,\delta[Y]$, it follows that $f_*[\mathcal{M}^1(X)]=\mathcal{M}^1(Y)$ so that $f_*$ is onto too. In fact, $f$ is continuous onto if and only if so is $f_*$.
In the case of $f\colon X\rightarrow X$ and $f_*\colon \mathcal{M}^1(X)\rightarrow \mathcal{M}^1(X)$, some dynamics, such as entropy and dimensions, of $f$ and $f_*$ may be described via each other (cf., e.g., \cite{BS75, GW95, BS25, SZ25, LL}).
\end{sse}

\subsection{Semi-open, almost 1-1 and irreducible maps}\label{s1.2}
\begin{sse}\label{1.2.1}
Let $\phi\colon X\rightarrow Y$ be a continuous onto mapping between topological spaces (not necessarily compact Hausdorff). We shall consider the following basic concepts: 
\begin{enumerate}[\textbf{a.}]
\item[\textbf{a.}] $\phi$ is called \textit{semi-open}~\cite{V70, B79}, or \textit{almost-open} \cite{Ak}, if $\mathrm{int}_Yf[U]\not=\emptyset$ for all $U\in\mathscr{O}(X)$; 
\item[\textbf{b.}] $\phi$ is called \textit{almost 1-1}, if 
$X_{\textrm{1-1}}[\phi]:=\{x\in X\,|\,\phi^{-1}(\phi(x))=\{x\}\}$ is dense in $X$ (cf., e.g.,~\cite{V70,A88});
\item[\textbf{c.}] $\phi$ is called \textit{irreducible}, if $A\in2^X$ with $\phi[A]=Y$ implies that $A=X$ (cf., e.g.,~\cite{E89, Wo, DeV}).
\end{enumerate}
Clearly, $\phi$ is irreducible if and only if for every $U\in\mathscr{O}(X)$ there is a point $y\in Y$ with $\phi^{-1}(y)\subseteq U$. So any almost 1-1 mapping is irreducible.
\end{sse}

If $X$ is a compact space and $Y$ a Hausdorff space here, then $\phi$ is irreducible if and only if 
for every $U\in\mathscr{O}(X)$ there exists $V\in\mathscr{O}(Y)$ with $\phi^{-1}[V]\subseteq U$ (cf., e.g.,~\cite{DeV}). Thus, an irreducible mapping is semi-open in that case. In particular, we have the following results:

\begin{sthm}\label{1.2.2}
Let $f\colon X\rightarrow Y$ be a continuous onto mapping between compact Hausdorff spaces. Then the following statements are satisfied: 
\begin{enumerate}[A.]
\item $f$ is irreducible if and only if so is $2^f$ (cf.~\cite[Thm.~9B]{DX}).

\item $f$ is irreducible if and only if so is $f_*$ (cf.~\cite[Thm.~1.2-($\mathcal{M}_i$)]{CH}).

\item $f$ is semi-open if and only if $2^f$ is semi-open (cf.~\cite[Thm.~4]{DX}).

\item If $f_*$ is semi-open, then so is $f$ (cf.~\cite[Thm.~B$^{\prime\prime}$]{DX}).
\end{enumerate}
\end{sthm}

\begin{sthm}\label{1.2.3}
Let $f\colon X\rightarrow Y$ be a continuous onto mapping between compact metric spaces. Then the following statements are satisfied:
\begin{enumerate}[A.]
    \item $f$ is almost 1-1 if and only if so is $2^f$ (cf.~\cite[Thm.~1.1-($\mathcal{H}_a^m$)]{CH}).
    \item $f$ is almost 1-1 if and only if so is $f_*$ (cf.~\cite[Thm.~1.2-($\mathcal{M}_a^m$)]{CH}).
    \item If $f$ is semi-open, then so is $f_*$ (cf.~\cite[Thm.~2.3]{G07}).
\end{enumerate}
\end{sthm}

``Semi-openness'' of continuous maps between compact Hausdorff spaces is essentially important for the structure theory of minimal dynamics (cf., e.g., \cite{V70, B79, Wo, A88, DeV, Ak}). In this paper we will be concerned with the following question based on Theorems~\ref{1.2.2} and \ref{1.2.3}:

\begin{sq}[cf.~\cite{DX, CH}]\label{1.2.4}
Let $f\colon X\rightarrow Y$ be a continuous onto map between compact Hausdorff non-metrizable spaces. Then:
\begin{enumerate}[a.]
    \item[\textbf{a.}] {\it If $f$ is semi-open, is $f_*\colon\mathcal{M}^1(X)\rightarrow\mathcal{M}^1(Y)$ semi-open?}
    \item[\textbf{b.}] {\it If $f_*\colon\mathcal{M}^1(X)\rightarrow\mathcal{M}^1(Y)$ is almost 1-1, is $f$ almost 1-1?}
\end{enumerate}
\end{sq}

Let $f\colon X\rightarrow Y$ be a continuous onto map between compact Hausdorff spaces. If $X$ is a metric space, then $Y$ is metrizable; but $X$ need not be metrizable if $Y$ is a metric space.
The following is one of our prime results proved in this paper, which is a weak solution to Question~\ref{1.2.4}a.

\begin{sthm}[Thm.~\ref{4.5}]\label{1.2.5}
    Let $f\colon X\rightarrow Y$ be a continuous onto map, where $X$ is a compact Hausdorff space and $Y$ a metric space. Then $f$ is semi-open if and only if so is $f_*$.
\end{sthm}

To prove Theorem~\ref{1.2.5}, we will present in $\S$\ref{s4} an inverse limit technique for describing the semi-openness of continuous mappings in the non-metrizable setting. We say that $f\colon X\rightarrow Y$ is ``quasi-separable'' if $f$ is the inverse limit of an inverse system $\{f_i\colon X_i\rightarrow Y\,|\,i\in\Lambda\}$ of continuous onto maps on a compact metric spaces (see Def.~\ref{4.1} for the details). We shall prove that if $f$ is quasi-separable, then $f$ is semi-open if and only if so is $f_*$.
In particular, if $Y$ is metrizable, then $f$ is quasi-separable (Lem.~\ref{4.2}) and so Question~\ref{1.2.4}a has a positive answer (Thm.~\ref{1.2.5}).

\subsection{Extensions of topological dynamics}\label{s1.3}
Beyond the setting of Theorem~\ref{1.2.5}, if $f$ preserves some dynamics, then we can show that $f$ and $f_*$ are both semi-open (Thm.~\ref{1.3.4}).

Let $S$ be a discrete monoid with identity element $e$ and $X$ a compact Hausdorff space. Then $S\curvearrowright X$ or $S\curvearrowright_\pi X$ is called a \textit{semiflow} (or an \textit{$S$-semiflow} if we need to emphasize the phase semigroup) \cite{EEN, AD}, denoted $\mathscr{X}$ if no confusion, if there exists a phase mapping $\pi\colon S\times X\rightarrow X$, $(t,x)\mapsto tx$, such that $ex=x$, $(st)x=s(tx)$ and the transit mapping $\pi_t\colon X\ni x\mapsto tx\in X$ is a continuous self-map of $X$, for all $s,t\in S$ and $x\in X$. Whenever $S$ is a group here, then $\mathscr{X}$ will be called a \textit{flow} (or an \textit{$S$-flow})~(cf., e.g.,~\cite{E69, B79, G76, Wo, A88, DeV}).
\begin{sse}\label{1.3.1}
Let $\mathscr{X}=S\curvearrowright_\pi X$ be a semiflow. Then, as usual (cf., e.g.,~\cite{B79, Wo, DeV, AD}), $\mathscr{X}$ is called: 
\begin{enumerate}[a.]
\item[\textbf{a.}] \textit{topologically transitive} (T.T.) if $X=\overline{S[U]}$ for all $U\in\mathscr{O}(X)$; 

\item[\textbf{b.}] \textit{minimal} if $\overline{Sx}=X$ for all $x\in X$; and a point $x\in X$ called \textit{almost periodic} (a.p.) for $\mathscr{X}$, if $\overline{Sx}$ is an $S$-minimal subset of $X$ (i.e., $S\curvearrowright\overline{Sx}$ is itself a minimal semiflow).
\end{enumerate}
In addition,
\begin{enumerate}[a.]
\item[\textbf{c.}] $\mathscr{X}$ is called a \textit{semi-open semiflow}, if $\pi_t\colon X\rightarrow X$ is a semi-open mapping for all $t\in S$. 
\end{enumerate}
\end{sse}

If $\mathscr{X}$ is minimal, then it is pointwise a.p.; an a.p. point is sometimes referred to as a uniformly recurrent point.  There is a topological criterion for the minimality of semiflows as follows:

\begin{1.3.1D}
A semiflow $S\curvearrowright X$ is minimal if and only if for every $\varepsilon\in\mathscr{U}_X$ there is a finite set $F\subseteq S$ such that $Fx$ is $\varepsilon$-dense in $X$.
\end{1.3.1D}

\noindent Indeed, sufficiency is obvious. Now, for necessity, let $\varepsilon, \alpha\in\mathscr{U}_X$ with $\alpha^3\subseteq\varepsilon$. Since $X$ is compact, there are finitely many points $x_1, \dotsc, x_n$ in $X$ with $\alpha[x_1]\cup\dotsm\cup\alpha[x_n]=X$. Moreover, as $x_i$ is a.p., it follows that there is a finite set $F$ in $S$ such that $Ftx_i\cap\alpha[x_i]\not=\emptyset$ for all $1\le i\le n$ and all $t\in S$. Then by $\overline{Sx_i}=X$ for $1\le i\le n$, we have that $\varepsilon[Fx]=X$ for all $x\in X$.

\begin{enumerate}[a.]
    \item[\textbf{e.}] Let $\mathscr{X}$ and $\mathscr{Y}$ be two $S$-semiflows. Then $\phi\colon\mathscr{X}\rightarrow \mathscr{Y}$ is referred to as an \textit{extension} of semiflows, if $\phi\colon X\rightarrow Y$ is a continuous onto mapping such that $\phi(tx)=t\phi(x)$ for all $x\in X$ and $t\in S$. Further, we say that $\mathscr{X}$ is 
\textit{$\phi$-fiber-onto} or $\phi$ is a \textit{fiber-onto extension} \cite[Def.~1.1.2c]{CD}, if $t\phi^{-1}(y)=\phi^{-1}(ty)$ for all $y\in Y$ and $t\in S$.
\end{enumerate}
Of course, an extension of flows is always fiber-onto. However, that of semiflows is generally not fiber-onto, even for a distal extension of $\mathbb{Z}_+$-semiflows (cf.~\cite[Ex.~1.1.3a]{CD}). 
 
Although a continuous onto mapping is generally not semi-open, an extension of minimal flows is always semi-open:

\begin{slem}[{cf.~\cite[Lem.~3.12.15]{B79} or \cite{V70, Wo, A88}}]\label{1.3.2}
If $\phi\colon \mathscr{X}\rightarrow \mathscr{Y}$ is an extension of flows where $\mathscr{Y}$ is minimal and $\mathscr{X}$ has a dense set of a.p. points, then $\phi\colon X\rightarrow Y$ is semi-open.
\end{slem}

The above simple observation is very useful for the structure theory of minimal topological dynamics. Even if $S\curvearrowright_\pi Y$ is a minimal semiflow, $\pi_t\colon y\in Y\mapsto ty\in Y$ need not be semi-open; and so, the standard proof of Lemma~\ref{1.3.2} (Thm.~\ref{P2.2.1}) is not valid even for extensions of minimal semiflows. However, using canonical commutative diagram (CD) of semiflows (Thm.~\ref{3.1}), in $\S$\ref{SO} we can extend Lemma~\ref{1.3.2} to semiflows as follows:

\begin{slem}[Thm.~\ref{3.2}]\label{1.3.3}
Let $\phi\colon \mathscr{X}\rightarrow \mathscr{Y}$ be a fiber-onto  extension of minimal semiflows. If $\mathscr{X}$ is metrizable, then $\phi\colon X\rightarrow Y$ and $\phi_*$ are semi-open mappings.
\end{slem}

In view of Glasner's Theorem~\ref{1.2.3}C, we essentially need to prove that $\phi$ is semi-open here.
Moreover, using canonical CD of flows (Thm.~\ref{3.1}) and Theorems~\ref{1.2.2} and \ref{1.3.2}, in $\S$\ref{SO} we can also extend Theorem~\ref{1.2.3}C to the non-metrizable setting as follows:

\begin{sthm}[Thm.~\ref{3.4}]\label{1.3.4}
If $\phi\colon \mathscr{X}\rightarrow \mathscr{Y}$ is an extension of minimal flows, then $\phi_*$ is semi-open.
\end{sthm}

\noindent Although $\pi_*\colon S\times\mathcal{M}^1(X)\rightarrow\mathcal{M}^1(X)$ is an affine flow and $\phi_*\colon\mathcal{M}^1(X)\rightarrow\mathcal{M}^1(Y)$ is an extension of affine flows in the setting of Theorem~\ref{1.3.4} (Rem.~\ref{5.1.6}), $S\curvearrowright\mathcal{M}^1(Y)$ is generally not minimal unless $Y$ is a singleton. Thus, Lemma~\ref{1.3.2} is not directly valid for the semi-openness of $\phi_*$ so that Theorem~\ref{1.3.4} is of interest itself.

The above two results (Lem.~\ref{1.3.3} $\&$ Thm.~\ref{1.3.4}) will be proved in in $\S$\ref{SO}. In fact, as a byproduct of proving them, we shall solve a question of Van der Woude \cite[Q.~IV.6.4b]{Wo} on the universal highly proximal extension of $\phi\colon \mathscr{X}\rightarrow \mathscr{Y}$ (Q.~\ref{3.7} $\&$ Cor.~\ref{3.9}).

Recall that $\phi\colon\mathscr{X}\rightarrow\mathscr{Y}$ is a ``quasi-separable extension'' iff $\phi$ is the inverse limit of an inverse system of extensions of compact metric semiflows (see Def.~\ref{5.2.1} for the details). Then $\phi$ is a quasi-separable extension in the setting of Lemma~\ref{1.3.3}. In view of that, Lemma~\ref{1.3.3} is contained in our another prime theorem that will be proved in $\S$\ref{s5.2}:

\begin{sthm}[Thm.~\ref{5.2.2}]\label{1.3.5}
 Let $\phi\colon\mathscr{X}\rightarrow\mathscr{Y}$ be a quasi-separable fiber-onto extension of minimal semiflows. Then $\phi$ and $\phi_*$ are both semi-open mappings.
\end{sthm}

In addition, in a quasi-separable extension $\phi\colon\mathscr{X}\rightarrow\mathscr{Y}$ with $Y=\{y\}$, $\mathscr{X}$ is called a quasi-separable semiflow.
In $\S$\ref{s5.3}, we shall deviate slightly from ``semi-openness'' and consider the structure of minimal weakly mixing quasi-separable flows (Thms.~\ref{Q5.10} $\&$ \ref{Q5.11}).

\subsection{Semi-openness of group homomorphisms}\label{s1.4}
Let $G$ be a group with a topology. $G$ is called a left-topological group, if $L_g\colon G\rightarrow G$, $x\mapsto gx$, is a continuous mapping for all $g\in G$. The following is a classic theorem of Pontryagin: 

\begin{sthm}[{Pontryagin's open-mapping theorem; cf.~\cite[Thm.~5.29]{HR63}}]\label{1.4.1}
  Let $G$ and $H$ be Hausdorff topological groups, where $G$ is locally compact Lindel\"{o}f and $H$ is Baire. If $f\colon G\rightarrow H$ is a continuous onto homomorphism, then $f$ is open.  
\end{sthm} 

\begin{sse}\label{1.4.2}
    Let $G$ and $H$ be two groups. A function $f\colon G\rightarrow H$ is said to be \textit{$G$-equivariant}, if there exists a function $\tau\colon G\rightarrow H$ such that $f(tx)=\tau(t)f(x)$ for all $t,x\in G$. In that case, it is easily verified that $\tau$ is a group homomorphism such that $f(x)=\tau(x)f(e)$ for all $x\in G$; and moreover, $\tau$ is surjective if so is $f$. If $f$ is itself a group homomorphism, then $f$ is $G$-equivariant with $\tau=f$.
\end{sse}

Finally, in $\S$\ref{s6} we shall generalize Pontryagin's Theorem using Lemma~\ref{1.3.2} (precisely, its generalization Thm.~\ref{P2.2.1}):

\begin{sthm}[Thm.~\ref{6.1}]\label{1.4.3}
Let $G$ and $H$ be left-topological groups, such that $G$ is Lindel\"{o}f quasi-regular (Def.~\ref{sP.2}b) and $H$ is Baire. If $f\colon G\rightarrow H$ is a locally closed (Def.~\ref{sP.2}c) continuous onto $G$-equivariant function, then $f$ is a semi-open mapping.
\end{sthm}

\begin{sthm}[Cor.~\ref{6.2}]\label{1.4.4}
Let $G$ and $H$ be left-topological groups, such that $G$ is locally compact Lindel\"{o}f pseudo-metrizable and $H$ is Hausdorff Baire. If $f\colon G\rightarrow H$ is a continuous onto $G$-equivariant function, then $f$ is an open mapping.
\end{sthm}
%%%%%%%%%%%%%%%%%%%%%%%%%%%%%%%%%%%%%%%%%%%%%%%%%%%%%%%%%%%%%%%%%%
\section{Preliminaries}\label{sP}%%%
In this section we shall introduce some basic notions and lemmas related to open, semi-open and irreducible mappings, which we shall need in our later and subsequent arguments.

\subsection{Open mappings}\label{sP.1}
Let $f\colon X\rightarrow Y$ be a map (continuous or not) between topological spaces. As usual, $f$ is called \textit{open}, if $f[U]\in\mathscr{O}(Y)$ for all $U\in\mathscr{O}(X)$.

Although our prime theorems are concerned with ``semi-openness'' of mappings, ``openness'' is sometimes useful for our dynamics approaches (e.g., in the later proof of Theorem~\ref{1.3.4} in $\S$\ref{SO}).
\begin{slem}\label{P2.1.1}
Let $f\colon X\rightarrow Y$ be any onto mapping between topological spaces. Then $f$ is open if and only if  $f[\mathbb{F}]$ is closed in $Y$ for every $F\in2^X$, where $\mathbb{F}=\{x\in X\,|\,f^{-1}(f(x))\subseteq F\}$.
\end{slem}

\begin{proof}
Indeed, for necessity, let $F\in 2^X$, $U=X\setminus F$; then $f[\mathbb{F}]=Y\setminus f[U]\in 2^Y$. Finally, for sufficiency, let $U\in\mathscr{O}(X)$ and $F=X\setminus U$; then  $f[U]=Y\setminus f[\mathbb{F}]$ is open so that $f$ is open.
\end{proof}

\begin{srem}\label{P2.1.2}
Lemma~\ref{P2.1.1} is a variant of Engelking \cite[Thm.~4.1.12]{E89} in which $f$ is assumed to be continuous. However, an open map is not necessarily continuous. For example, let $X$ be a set equipped with topologies $\mathfrak{T}_1$ and $\mathfrak{T}_2$ with $\mathfrak{T}_1\varsubsetneq\mathfrak{T}_2$; then $\textrm{id}_X\colon(X,\mathfrak{T}_1)\rightarrow(X,\mathfrak{T}_2)$ is an open, closed, 1-1, and onto map, but it is not continuous.
\end{srem}

\begin{sse}\label{P2.1.3}
Let $f\colon X\rightarrow Y$ be an onto map between topological spaces. Then its \textit{adjoint mapping} $f_\textsl{ad}\colon Y\rightarrow2^X$, defined by $y\mapsto \overline{f^{-1}(y)}$, is called \textit{lower semicontinuous}, if for every net $\{y_\alpha\,|\,\alpha\in A\}$ with $y_\alpha\to y$ in $Y$ we have that $f^{-1}(y)\subseteq\bigcap_{\alpha\in A}\overline{\bigcup\{f^{-1}(y_i)\colon i\ge\alpha\}}$. 
\end{sse}

\begin{slem}[{cf.~\cite[Thm.~II.1.3d]{Wo} for $f$ a continuous onto map of compact Hausdorff spaces}]\label{2.1.4}
$f$ is open if and only if $f_\textsl{ad}\colon Y\rightarrow 2^X$ is lower semicontinuous.
\end{slem}

\begin{proof}
Let $f$ be open and $y_\alpha\to y$ in $Y$. Then by Lemma~\ref{P2.1.1}, $F_\alpha:=\overline{\bigcup\{f^{-1}(y_i)\colon i\ge \alpha\}}\in2^X$ such that $y\in f[\mathbb{F}_\alpha]$; and so, $f^{-1}(y)\subseteq F_\alpha$ for all $\alpha\in A$ and $f^{-1}(y)\subseteq\bigcap_{\alpha}F_\alpha$. Conversely, to the contrary, let $x\in X$ and $U\in\mathfrak{N}_x(X)$ such that $f[U]\notin\mathfrak{N}_{f(x)}(Y)$. Then there is a net $y_\alpha\notin f[U]\to f(x)$ in $Y$; and so, $f^{-1}(y_\alpha)\subseteq X\setminus U$ for all $\alpha$. Thus, by sufficiency condition, $x\in f^{-1}(f(x))\subseteq X\setminus U$, contrary to $x\in U$. The proof is complete.
\end{proof}

Note that if $f\colon X\rightarrow Y$ is a continuous map of a compact space $X$ onto a Hausdorff space $Y$, then $f_{ad}\colon Y\rightarrow 2^X$ is always upper semi-continuous (i.e., for all $y\in Y$ and $U\in\mathfrak{N}_{f^{-1}(y)}(X)$, there exists a set $V\in\mathfrak{N}_y(Y)$ such that $f^{-1}[V]\subseteq U$).

Consequently, if $f\colon X\rightarrow Y$ is a continuous onto map between compact Hausdorff spaces, then by Lemmas~\ref{P2.1.1} and \ref{2.1.4}, $f$ is open, if and only if $f_\textsl{ad}\colon Y\rightarrow2^X$ is continuous, if and only if a net $y_n\to y$ in $Y$ implies that for all $x\in f^{-1}(y)$ there exists a net $x_n\in f^{-1}(y_n)\to x$ in $X$.

The ``openness'' of a continuous onto map between compact Hausdorff spaces may then be characterized via its induced mappings as follows:

\begin{sthm}\label{2.1.5}
Let $f\colon X\rightarrow Y$ be a continuous onto map between compact Hausdorff spaces; then the following two statements are satisfied:
\begin{enumerate}[(1)]
\item $f$ is open if and only if $f_*$ is open (cf.~\cite[Prop.~4.1 $\&$ Thm.~4.4]{DE}).
\item $f$ is open if and only if $2^f$ is open (\cite[Thm.~3]{DX}; cf.~\cite[Thm.~4.3]{H97} for $X, Y$ in the category of continua).
\end{enumerate}
\end{sthm}

\subsection{Semi-open mappings}\label{sP.2}

Let $\mathscr{X}$ be an $S$-flow, where $X$ need not be a compact Hausdorff space. Recall that $\mathscr{X}$ is minimal iff $\overline{Sx}=X$ for all $x\in X$; and, a point $x\in X$ is a.p. iff $\overline{Sx}$ is an $S$-minimal subset of $X$ (cf.~Def.~\ref{1.3.1}b). 

\begin{enumerate}[\ref{sP.2}a]
\item $X$ is a Lindel\"{o}f space \cite{K55, E89} iff every open cover has a countable subcover. 

\item $X$ is quasi-regular \cite{O60, M75} iff for every $U\in\mathscr{O}(X)$ there is a set $V\in\mathscr{O}(X)$ with $\bar{V}\subseteq U$. 

\item A function $f\colon X\rightarrow Y$ between topological spaces is called \textit{locally closed}, if for all $U\in\mathscr{O}(X)$ there exists a set $V\in\mathscr{O}(U)$ such that $f|_V\colon V\rightarrow Y$ is a closed mapping.
\end{enumerate}

In the class of locally compact spaces, a space has the Lindel\"{o}f property if and only if it is $\sigma$-compact. 
Next we can simply generalize and prove Lemma~\ref{1.3.2} as follows:

\begin{sthm}\label{P2.2.1}
Let $\phi\colon\mathscr{X}\rightarrow\mathscr{Y}$ be a locally closed extension of $S$-flows, where $X$ is a quasi-regular Lindel\"{o}f space and $Y$ a Baire space ($X,Y$ are not necessarily locally compact).
If $\mathscr{X}$ has a dense set of a.p. points and $\mathscr{Y}$ is minimal, then $\phi\colon X\rightarrow Y$ is semi-open.
\end{sthm}

\begin{proof}
    Since $\phi[\overline{Sx}]=Y$ for all $x\in X$, we can assume that $\mathscr{X}$ is minimal without loss of generality. Let $U\in\mathscr{O}(X)$. Since $X$ is quasi-regular and $\phi$ is locally closed, we can choose some set $V\in\mathscr{O}(X)$ with $\bar{V}\subseteq U$ such that $\phi[\bar{V}]$ is closed in $Y$. By minimality and Lindel\"{o}f property, there exists a countable set $\{t_n\colon n\in\mathbb{N}\}$ in $S$ such that $\bigcup_nt_n^{-1}\bar{V}=X$. Then $\textrm{int}_Y\phi[t_i^{-1}\bar{V}]\not=\emptyset$ for some $i\in\mathbb{N}$. Note that $y\in Y\mapsto t_iy\in Y$ is a homeomorphism and so open, since $S\curvearrowright Y$ is a flow. As $t_i\phi[t_i^{-1}\bar{V}]=\phi[\bar{V}]\subseteq\phi[U]$, it follows that $\textrm{int}_Y\phi[U]\not=\emptyset$. Thus, $\phi\colon X\rightarrow Y$ is semi-open. 
\end{proof}

Note that if $f\colon X\rightarrow Y$ is a continuous map of a locally compact quasi-regular space $X$ into a Hausdorff space,
then $f$ is locally closed. In view of that, Theorem~\ref{P2.2.1} conains Lemma~\ref{1.3.2}.

\begin{sse}[Semi-open sets and semicontinuous/quasicontinuous mappings]\label{2.2.2}
Let $f\colon X\rightarrow Z$ be a mapping (continuous or not) between topological spaces. 
%By $D(f)$ we denote the set of points of discontinuity of $f$.
\begin{enumerate}[a.]
\item[\textbf{a.}] A set $A$ in a space $X$ is termed \textit{semi-open} \cite[Def.~1]{L63}, if there exists a set $U\in\mathscr{O}(X)$ with $U\subseteq A\subseteq\bar{U}$. Then a subset $A$ of $X$ is semi-open if and only if $A\subseteq\overline{\mathrm{int}\,A}$ 
(\cite[Thm.~1]{L63}).

\item[\textbf{b.}] Following \cite{K32, B52, HT}, $f$ is called \textit{quasicontinuous at a point $x\in X$}, if to each $U\in\mathfrak{N}_{x}(X)$ and $V\in\mathfrak{N}_{f(x)}(Z)$, there exists $W\in\mathscr{O}(U)$ such that $f[W]\subseteq V$. We say that $f$ is \textit{quasicontinuous}, if $f$ is quasicontinuous at every point of $X$.

\item[\textbf{c.}] Following \cite[Def.~4]{L63}, $f$ is termed \textit{semicontinuous}, if $f^{-1}[V]$ is semi-open in $X$ for every $V\in\mathscr{O}(Z)$.
\end{enumerate}

\begin{2.2.2D}
Let $f\colon X\rightarrow Z$ be a map between topological spaces. Then $f$ is semicontinuous if and only if it is quasicontinuous.
\end{2.2.2D}

\begin{proof}
    It is obvious by Defs.~\ref{2.2.2}b and \ref{2.2.2}c. 
\end{proof}

For example, if $L$ is the Sorgenfrey line (cf.~\cite[$\S$5.3, Ex.~1]{W70}) and $\mathbb{R}$ the 1-dimensional euclidean space, then $\textrm{id}\colon\mathbb{R}\rightarrow L$ is an open noncontinuous and semicontinuous mapping. The following criterion is motivated by Lemma~\ref{2.1.4}.

\begin{2.2.2E}
Let $f\colon X\rightarrow Z$ be a continuous onto mapping between compact Hausdorff spaces. If $f_\textsl{ad}\colon Z\rightarrow 2^X$ is semicontinuous, then $f$ is semi-open.
\end{2.2.2E}

\begin{proof}
%First we note that since $f\colon X\rightarrow Z$ is a continuous %onto mapping between compact Hausdorff spaces, hence %$f_\textsl{ad}\colon Z\rightarrow 2^X$ is upper semicontinuous.
Let $U\in\mathscr{O}(X)$ be arbitrarily given. If $U$ contains a fiber of $f$, then $\textrm{int}f[U]\not=\emptyset$. Next, suppose that $U$ does not contain any fiber of $f$. To prove that $\textrm{int}f[U]\not=\emptyset$, take a set $V\in\mathscr{O}(X)$ such that $\bar{V}\subseteq U$. Then $\{U,X\setminus\bar{V}\}$ is an open cover of $X$ such that $f^{-1}(z)\in\langle U,X\setminus\bar{V}\rangle$ for all $z\in f[U]$. If $\textrm{int}f[U]=\emptyset$, then for all $z_0\in f[U]$, there exists (and if there is) a net $z_\alpha\notin f[U]$ with $z_\alpha\to z_0$; and so (then), $f^{-1}(z_0)\in\langle U,X\setminus\bar{V}\rangle$ and $f^{-1}(z_\alpha)\cap U=\emptyset$ for all $\alpha$; this implies that $f^{-1}(z_\alpha)\notin\langle U,X\setminus\bar{V}\rangle$ for all $\alpha$, contrary to the hypothesis, this is because $f[U]$ is contains in $\textrm{cl}_Z\textrm{int}_Zf_\textsl{ad}^{-1}[\langle U,X\setminus\bar{V}\rangle]$. The proof is complete.
\end{proof}

\begin{2.2.2F}[Semi-openness of $f$ $\not\Rightarrow$ semicontinuity of $f_\textsl{ad}$]
Let $X=Z\times\{1\}\sqcup Z\times\{2\}$ be Ellis' two-circle space \cite[Ex.~5.29]{E69}, where $Z=\mathbb{S}$ is the unit circle in $\mathbb{C}$ oriented in a counter-clockwise direction. Then for all $z\in Z$, we have that 
\begin{enumerate}
    \item[] $\mathfrak{N}_{(z,1)}(X)=\{[z,a)\times\{1\}\cup(z,a]\times\{2\}\,|\,a\in\mathbb{S}\setminus\{z\}\}$,
    \item[] $\mathfrak{N}_{(z,2)}(X)=\{[b,z)\times\{1\}\cup(b,z]\times\{2\}\,|\,b\in\mathbb{S}\setminus\{z\}\}$.
\end{enumerate}
Let $f\colon X\rightarrow Z$ be the canonical projection. Clearly, $f$ is semi-open. However, $f_\textsl{ad}\colon Z\rightarrow 2^X$ is not semicontinuous. Indeed, let $z\in Z$; then for all $a,b\in\mathbb{S}\setminus\{z\}$ such that the lengths of the arcs $[b,z]$ and $[z,a]$ are less than 1/2, we have that 
$$
f_\textsl{ad}(z)=f^{-1}(z)=\{(z,1), (z,2)\}\in\langle U_1,U_2\rangle\ \left(=\left\{A\in2^X\,|\,A\subseteq U_1\cup U_2\ \&\ A\cap U_1\not=\emptyset\not=A\cap U_2\right\}\right),
$$
where $U_1=[z,a)\times\{1\}\cup(z,a]\times\{2\}$ and $U_2=[b,z)\times\{1\}\cup(b,z]\times\{2\}$; but, for all $z^\prime\not=z$ in $Z$, $f_\textsl{ad}(z^\prime)=f^{-1}(z^\prime)\notin\langle U_1,U_2\rangle$. Thus, $f_\textsl{ad}$ is not quasicontinuous (semicontinuous).
\end{2.2.2F}

Clearly, ``semi-open mapping'' and ``semicontinuous mapping'' are two different concepts. However, if $f$ is a 1-1 onto mapping, then $f$ is semi-open if and only if $f^{-1}$ is semicontinuous. This simple observation is useful for topological groups. For example, if $G$ is a group on a topological space and $\Box^{-1}\colon G\ni x\mapsto x^{-1}\in G$, then $\Box^{-1}$ is semi-open if and only if $\Box^{-1}$ is semicontinuous.
\end{sse}

\begin{slem}\label{P2.2.3}
Let $f\colon X\rightarrow Y$ be a continuous mapping between topological spaces. Then: 
\begin{enumerate}[(1)]
    \item $f$ is semi-open if and only if the preimage of every dense subset of $Y$ is dense in $X$ (cf.~\cite[Lem.~2.1]{G18}; in fact, the continuity of $f$ plays no role here).

    \item $f$ is semi-open if and only if $\mathrm{int}_Yf[U]$ is dense in $f[U]$ for all $U\in\mathscr{O}(X)$.
\end{enumerate}
\end{slem}

\begin{proof}
We need only prove (2). For that, sufficiency is obvious. Now suppose $f$ is semi-open. Let $U\subseteq X$ be an open set. As $U\setminus f^{-1}\left[\overline{\mathrm{int}\,f[U]}\right]$ is open in $X$, it follows that $U\setminus f^{-1}\left[\overline{\mathrm{int}\,f[U]}\right]=\emptyset$ so that $f[U]\subseteq\overline{\mathrm{int}\,f[U]}$. Thus, $\mathrm{int}\,f[U]$ is dense in $f(U)$. The proof is complete.
\end{proof}

Thus, a continuous map $f\colon X\rightarrow Y$ is semi-open, iff $f[U]$ is semi-open for all $U\in\mathscr{O}(X)$, iff the image of every semi-open set is semi-open. Here the continuity of $f$ has played a role.

\begin{sse}\label{2.2.4}
A family $\mathcal{B}\subset\mathscr{O}(X)$ is a \textit{$\pi$-base} for $X$ \cite{O60, W70, M75}, if any $U\in\mathscr{O}(X)$ contains some member of $\mathcal{B}$.
Let $p\colon X\rightarrow Y$ be a continuous mapping; then we say that $X$ has a \textit{$p$-fiber countable $\pi$-base}, if there exists a countable family $\mathcal{U}\subseteq\mathscr{O}(X)$ such that $\{U\cap p^{-1}(y)\,|\,U\in\mathcal{U}\}$ is a $\pi$-base for $p^{-1}(y)$
for all $y\in Y$.
\end{sse}
A second countable space has of course a $\pi$-base; but not vice versa. For example, the Stone-\v{C}ech compactification 
$\beta\mathbb{N}$ of $\mathbb{N}$ is not a second countable space, but it has a countable $\pi$-base, say $\mathcal{B}=\{\{b\}\colon b\in\mathbb{N}\}$.

As an application of ``semi-openness'' condition, we shall prove a topological Fubini theorem as follows:

\begin{sthm}[{cf.~\cite[Lem.~5.3]{D23} for $X$ a Polish space}]\label{P2.2.5}
Suppose that $p\colon X\rightarrow Y$ is a semi-open continuous onto mapping, where
$X$ has a $p$-fiber countable $\pi$-base $\mathscr{U}$. If $G\subseteq X$ is a dense open set, then 
$Y_G=\left\{y\in Y\,|\,G\cap p^{-1}(y)\textrm{ is dense open in }p^{-1}(y)\right\}$
is residual in $Y$. In particular, if $K\subseteq X$ is residual, then
$Y_K=\left\{y\in Y\,|\,K\cap p^{-1}(y)\textrm{ is residual in }p^{-1}(y)\right\}$
is residual in $Y$.
\end{sthm}

\begin{proof}
Let $F=X\setminus G$. Then $F$ is a nowhere dense closed set in $X$. Let $F_y=F\cap p^{-1}(y)$ for all $y\in Y$. Let
$B=\left\{y\in Y\,|\,\mathrm{int}_{p^{-1}(y)}F_y\not=\emptyset\right\}$.
So if $y\notin B$, then $G_y$ is open dense in $p^{-1}(y)$. Thus, $Y\setminus B\subseteq Y_G$ and we need only prove that $B$ is meager in $Y$.
For that, write $\mathscr{U}=\{U_n\}_{n=1}^\infty$. If $y\in B$, then $U_n\cap p^{-1}(y)\subseteq F_y$ for some $n\in\mathbb{N}$.
Put
$C_n=\{y\in B\,|\,U_n\cap p^{-1}(y)\subseteq F_y\}$ and $D_n=\mathrm{int}_Y \bar{C}_n$
for all $n\in\mathbb{N}$. Then $B=\bigcup_{n=1}^\infty C_{n}$, and $B$ is meager in $Y$ if each $D_{n}=\emptyset$. Indeed, if $D_{n}\not=\emptyset$, then $U_n\cap p^{-1}(y)\subseteq F_y$ for all $y\in D_n\cap C_n$ and $D_n\cap C_n$ is dense in $D_n$. So $U_n\cap p^{-1}[D_n\cap C_n]\subseteq F$. Further, by Lemma~\ref{P2.2.3}, it follows that $\emptyset\not=U_n\cap p^{-1}[D_n]\subseteq \bar{F}=F$, contrary to $F$ being nowhere dense in $X$. The proof is complete.
\end{proof}

We notice that the above Fubini theorem is due to L.\,E.\,J. Brouwer 1919 for the special case $p\colon [0,1]\times[0,1]\xrightarrow{(x,y)\mapsto x}[0,1]$, to C.~Kuratowski and S.~Ulam 1932 for $p\colon X\times Y\xrightarrow{(x,y)\mapsto x}X$ where $X,Y$ are separable metric spaces, and to Oxtoby 1960 \cite{O60} for $p\colon X\times Y\xrightarrow{(x,y)\mapsto x}X$ with $Y$ having a countable $\pi$-base. See Veech (1970) \cite[Prop.~3.1]{V70} and Glasner (1990) \cite[Lem.~5.2]{G90} for the case that $p$ is an extension of a minimal compact metric flow $\mathscr{X}$ to a minimal flow $\mathscr{Y}$.

\begin{sse}
Recall that a sequence $\{\mathscr{U}_n\}_{n=1}^\infty$ of open covers of a space $X$ is a \textit{development} for $X$, if the family $\{\textrm{st}(x,\mathscr{U}_n)=\cup\{U\in\mathscr{U}_n\,|\,x\in U\}\colon n\in\mathbb{N}\}$ is a base for $\mathfrak{N}_x(X)$, for each $x\in X$. A \textit{developable space} is a space which has a development; and a regular developable space is called a \textit{Moore space} (cf.~\cite[Def.~1.3]{G84}). This is a generalization of pseudo-metric spaces.
\end{sse}

\begin{slem}\label{P2.2.7}
Let $p\colon X\rightarrow Y$ be a semi-open continuous onto mapping, where $X$ is a developable space. Then
$X_o[p]=\{x\in X\,|\,p[U]\in\mathfrak{N}_{p(x)}(Y)\ \forall U\in\mathfrak{N}_x(X)\}$ is a residual set in $X$.
\end{slem}

\begin{proof}
Given any $n\in\mathbb{N}$, let $X_n=\{x\in X\,|\,\exists U\in\mathscr{U}_n\textrm{ s.t. } U\in\mathfrak{N}_x(X)\ \&\ p[U]\in\mathfrak{N}_{p(x)}(Y)\}$. Then by Lemma~\ref{P2.2.3}, it follows that $X_n$ is open dense in $X$. Thus, $X_o[p]=\bigcap_{n=1}^\infty X_n$ is a residual set in $X$. The proof is complete.
\end{proof}

Lemma~\ref{P2.2.7} will be needed in $\S$\ref{s6}. In addition, we shall need the following well-known result in Theorem~\ref{3.1} and Example~\ref{4.8}, for which we give a simple proof without using the $\varepsilon$-spanning or $\varepsilon$-separating numbers of $\psi(y)$ here.

\begin{slem}[{cf.~\cite[Thm.~II.1.3e]{Wo} or \cite[Lem.~14.44]{A88}}]\label{P2.2.8}
Let $\phi\colon X\rightarrow Y$ be a continuous onto map between compact metric spaces. Then the set $D(\phi_{\textsl{ad}})$ of points of discontinuity of $\phi_\textsl{ad}$
is meager in $Y$. (Thus, $\phi$ is irreducible if and only if $\phi$ is almost 1-1 in the category of compact metric spaces.)
\end{slem}

\begin{proof}
Let $Z=\phi_{\textsl{ad}}[Y]\subset2^X$. Then $\phi_\textsl{ad}\colon Y\rightarrow Z$ is 1-1 open ($\because\forall V\in\mathscr{O}(Y)$, $\phi_{\textsl{ad}}[V]=Z\cap\langle \phi^{-1}[V]\rangle$); and $\phi_{\textsl{ad}}^{-1}\colon Z\subset2^X\rightarrow Y$, $\phi^{-1}(y)\mapsto y$ is uniformly continuous. So, $\phi_{\textsl{ad}}^{-1}$ admits a unique uniformly continuous extension, denoted $\widetilde{\phi_\textsl{ad}^{-1}}\colon \bar{Z}\rightarrow Y$\,(cf.~\cite[Thm.~6.26]{K55}).

To prove that $Y_c(\phi_\textsl{ad})=Y\setminus D(\phi_\textsl{ad})$ is residual in $Y$, it suffices to prove that for all $\varepsilon>0$, $Y_{\varepsilon,c}(\phi_\textsl{ad})$---the set of $\varepsilon$-continuous points of $\phi_\textsl{ad}$ is dense, since $Y_{\varepsilon,c}(\phi_\textsl{ad})$ is always open. For that, let $V$ and $U$ be open sets in $Y$ with $\emptyset\not=V\subseteq\bar{V}\subseteq U$. Since $\phi_\textsl{ad}[V]$ is open in $Z$, there is a sequence of open balls $\{B_n\}_{n=1}^\infty$ in $Z$ with diameter
$|B_n|<\varepsilon/2$ and $\phi_\textsl{ad}[V]=\bigcup_{n=1}^\infty B_n$. Then
$$
V={\bigcup}_{n=1}^\infty \phi_\textsl{ad}^{-1}[B_n]={\bigcup}_{n=1}^\infty \widetilde{\phi_\textsl{ad}^{-1}}[B_n]\subseteq {\bigcup}_{n=1}^\infty \overline{\widetilde{\phi_\textsl{ad}^{-1}}[B_n]}={\bigcup}_{n=1}^\infty \widetilde{\phi_\textsl{ad}^{-1}}[\bar{B}_n]\subseteq\bar{V},
$$
where $\bar{B}_n$ is the closure of $B_n$ in $\bar{Z}$. Let $W_n=\mathrm{int}\,\widetilde{\phi_\textsl{ad}^{-1}}[\bar{B}_n]$. By Baire's theorem, $W_n\not=\emptyset$ for some $n\in\mathbb{N}$. Thus, $\phi_\textsl{ad}[W_n]\subseteq\bar{B}_n$ and $|\phi_\textsl{ad}[W_n]|<\varepsilon$. This shows that $Y_{\varepsilon,c}(\phi_\textsl{ad})\cap U\not=\emptyset$. Since $U$ is arbitrary, $Y_{\varepsilon,c}(\phi_\textsl{ad})$ is dense in $Y$. The proof is complete.
\end{proof}

Note that if $\phi$ is semi-open (such as in Lemma~\ref{1.3.2} and Theorem~\ref{1.3.3}), then Lemma~\ref{P2.2.8} follows easily from Lemma~\ref{P2.2.7} and Theorem~\ref{P2.2.5}.

\begin{srem}
Using a proof same as that of \cite[Thms.~1 $\&$~2]{Fo51}, we can improve Lemma~\ref{P2.2.8} as follows:

\begin{xthm}[{cf.~\cite[Thms.~1 $\&$~2]{Fo51} for $X$ a metric space}]
Let $\psi\colon Y\rightarrow 2^X$ be either upper semi-continuous or lower semi-continuous, where $X$ is a pseudo-metric space. Then the set
$D(\psi)$ of points of discontinuity of $\psi$ is meager in $Y$.
\end{xthm}
\end{srem}

\subsection{Highly proximal extensions}\label{sP.3}%%%%%
The following technical concepts and lemmas are mainly for proving our prime theorem---Theorem~\ref{1.3.4} in $\S$\ref{SO}. 

\begin{sse}[$\diamond$ operation]\label{P2.3.1}
Let $\phi\colon\mathscr{X}\rightarrow\mathscr{Z}$ be an extension of $S$-semiflows. By $\beta S$ it means the Stone-\v{C}ech compactification of $S$ with a right-topological semigroup structure.

 The ``circle operation'' (cf.~\cite{G76, Wo, A88, DeV} or \cite[A.1.1]{CD}), as the extension to $\beta S$ of the hyperspace semiflow $2^\mathscr{X}:=S\curvearrowright2^X$, is defined as follows:
\begin{enumerate}
\item[] $p\diamond K=\lim_i t_iK\in 2^X$ for all $p\in\beta S$ and $K\in2^X$, where $t_i\in S\to p$ in $\beta S$.
\end{enumerate}
Here $p\diamond K$, as a ``point'' of $2^X$, is independent of the choice of the net $\{t_i\}$ in $S$ satisfying $t_i\to p$ in $\beta S$. 
\begin{xnote}
  In the literature, $p\diamond K$ has been written as $p\circ K$ or $p^\circ K$ or $p\odot K$ and so on. Here using symbol $\diamond$ is to distinguish from the composition of mappings.  
\end{xnote}

Next we can define a closed subspaces of the hyperspace $2^X$ associated with $\phi\colon\mathscr{X}\rightarrow\mathscr{Z}$ as follows:
\begin{enumerate}
\item[] $2^{X,\phi}=\{A\in 2^X\,|\,\exists z\in Z\textrm{ s.t. }A\subseteq\phi^{-1}(z)\}$.
\end{enumerate}
Clearly, $2^{X,\phi}$ is an $S$-invariant subset of $2^X$ so that
\begin{enumerate}
\item[] $2^{\mathscr{X},\phi}:=S\curvearrowright 2^{X,\phi}$
\end{enumerate}
is a subsemiflow of $2^\mathscr{X}$. Let
\begin{enumerate}
\item[] $\tilde{\chi}\colon 2^{X,\phi}\rightarrow Z$ be defined by $A\in2^{X,\phi}\mapsto z\in Z$ iff $A\subseteq\phi^{-1}(z)$.
\end{enumerate}
Clearly, $\tilde{\chi}$ is $S$-equivariant so that $\mathscr{Z}$ is a factor of $2^{\mathscr{X},\phi}$ by means of $\tilde{\chi}$. In particular, we have for $z\in Z$ and $p\in \beta S$ that $\tilde{\chi}(p\diamond\phi^{-1}(z))=pz$.
\end{sse}

\begin{sse}[Highly proximal extension]\label{P2.3.2}
Let $\phi\colon\mathscr{X}\rightarrow\mathscr{Z}$ be an extension of $S$-semiflows, where $\mathscr{X}$ and $\mathscr{Z}$ need not be minimal.
As in \cite{E73, S76, AG77, Wo, DeV}, 
\begin{enumerate}[a.]
\item[\textbf{a.}] $\phi$ is \textit{proximal} or $\mathscr{X}$ is a \textit{proximal extension} of $\mathscr{Z}$ via $\phi$, iff $\varDelta_X\cap \overline{S(x_1,x_2)}\not=\emptyset$ $\forall x_1,x_2\in\phi^{-1}(z)$ for all $z\in Z$, where $\varDelta_X$ is the diagonal of $X\times X$;

\item[\textbf{b.}] $\phi$ is \textit{highly proximal} (h.p.) or $\mathscr{X}$ is an \textit{h.p. extension} of $\mathscr{Z}$ via $\phi$, iff for all $z\in Z$ there is a net $t_n\in T$ such that $t_n[\phi^{-1}(z)]$ converges to a singleton in $2^X$, iff for all $z\in Z$ there is $p\in \beta S$ with $p\diamond\phi^{-1}(z)=\{px\}$ for all $x\in\phi^{-1}(z)$. In this case, $\phi$ is of course proximal.
\end{enumerate}
\end{sse}

Let $\mathscr{Z}$ be a minimal semiflow. If $\phi$ is almost 1-1, then it is obviously an h.p. extension. If $\mathscr{X}$ is a $\phi$-fiber-onto metrizable semiflow, then $\phi$ is h.p. if and only if it is almost 1-1 by Lemmas~\ref{P2.2.8} and \ref{P2.3.3} below (see, e.g., \cite[Rem.~IV.1.2]{Wo} for $S$ a group).

\begin{slem}[{cf.~\cite{AG77}, \cite[Thm.~IV.2.3]{Wo}, \cite[Prop.~VI.3.3]{DeV} for minimal flows}]\label{P2.3.3}
Let $\phi\colon\mathscr{X}\rightarrow\mathscr{Z}$ be a fiber-onto extension of semiflows, where $\mathscr{Z}$ is minimal and $\mathscr{X}$ has a dense set of a.p. points. Then the following conditions are pairwise equivalent:
\begin{enumerate}[(1)]
\item $\phi$ is h.p. (so $\mathscr{X}$ is minimal);
\item $\phi$ is irreducible;
\item $\tilde{\chi}\colon2^{\mathscr{X},\phi}\rightarrow\mathscr{Z}$ is h.p.;
\item $\tilde{\chi}\colon2^{\mathscr{X},\phi}\rightarrow\mathscr{Z}$ is proximal;
\item $2^{\mathscr{X},\phi}$ has a unique $S$-minimal subset.
\end{enumerate}
\end{slem}

\begin{proof}
$(1)\Rightarrow(2)$: Let $U\in\mathscr{O}(X)$ and $x_0\in U$ an a.p. point for $\mathscr{X}$. Since $\phi$ is h.p., there exists an element $p\in\beta S$ with $p\diamond\phi^{-1}(z_0)=\{px_0\}$, where $z_0=\phi(x_0)$. As $x_0$ is a.p., it follows that there is some element $q\in\beta S$ with $q\diamond\phi^{-1}(z_0)=\{x_0\}$ so that $U$ contains a fiber of $\phi$. Thus, $\phi$ is irreducible by Def.~\ref{1.2.1}c. (``$\phi$-fiber-onto'' plays a role only here.)

$(2)\Rightarrow(3)$: Let $z_0\in Z$ and $x_0\in\phi^{-1}(z_0)$. Since $\mathscr{Z}$ is minimal and $\phi$ is irreducible, there exists a net $t_n\in S$ such that $t_n\phi^{-1}(z_0)\to\{x_0\}$. However, every $A\in\tilde{\chi}^{-1}(z_0)$ is contained in $\phi^{-1}(z_0)$. Hence $t_n\tilde{\chi}^{-1}(z_0)\to\{x_0\}$. Thus, $\tilde{\chi}$ is h.p.

$(3)\Rightarrow(4)$: Obvious.

$(4)\Rightarrow(5)$: Obvious. 

$(5)\Rightarrow(1)$: Clearly, $\{\{x\}\,|\,x\in X\}$ is an $S$-invariant closed subset of $2^{\mathscr{X},\phi}$. Let $z\in Z$ and $\mathfrak{z}=\phi^{-1}(z)\in 2^{X,\phi}$. As $\overline{S\mathfrak{z}}$ contains an a.p. point of $2^{\mathscr{X},\phi}$, it follows that there exists some element $p\in\beta S$ such that $p\diamond\phi^{-1}(z)$ is a singleton. Thus, $\phi$ is h.p. The proof is complete.
\end{proof}

\begin{slem}[{Ellis-Shoenfeld-Auslander-Glasner; cf.~\cite{S76, AG77, Wo} in the class of minimal flows}]\label{P2.3.4}
Let $\phi\colon\mathscr{X}\rightarrow\mathscr{Z}$ be an extension of semiflows with $\mathscr{Z}$ being minimal. Let
\begin{enumerate}
\item[] $Z_\phi=\mathrm{cl}_{2^{X,\phi}}\{t\phi^{-1}(z)\,|\,z\in Z\ \&\  t\in S\}\subseteq2^{X,\phi}$.
\end{enumerate}
Then:
\begin{enumerate}[(1)]
\item $\chi=\tilde{\chi}|_{Z_\phi}\colon Z_\phi\rightarrow Z$ is h.p.;
\item $Z_\phi$ has a unique $S$-minimal subset, denoted $Z_\phi^\natural$;
\item If $\mathscr{X}$ is minimal, then $\phi$ is h.p. if and only if $\mathscr{X}\rightarrow\mathscr{Z}_\phi^\natural$, $x\mapsto\{x\}$ is an isomorphism;
\item If $\mathscr{X}$ has a dense set of a.p. points, then $\phi$ is fiber-onto open if and only if $\chi\colon\mathscr{Z}_\phi^\natural\rightarrow\mathscr{Z}$ is an isomorphism.
\end{enumerate}
If $X$ is a metric space, then so is $Z_\phi^\natural$.
\begin{xnote}
$\chi\colon \mathscr{Z}_\phi^\natural\rightarrow\mathscr{Z}$  is called the h.p. quasifactor representation of $\phi\colon\mathscr{X}\rightarrow\mathscr{Z}$ in $\mathscr{X}$.
\end{xnote}
\end{slem}

\begin{proof}
    First of all, we note that $Z_\phi$ is an $S$-invariant closed subset of $2^{\mathscr{X},\phi}$. By Zorn's Lemma, we can select a minimal point in $Z_\phi$ under inclusion. 
    
(1): Let $A\in Z_\phi$ be a minimal point. Then we can find a point $z\in Z$ and an element $p\in\beta S$ with 
$p\diamond\phi^{-1}(z)=A$ and a net $t_i\in S$ with $t_i\to p$ in $\beta S$. This implies that $t_i\chi^{-1}(z)\to\{A\}$ in $2^{Z_\phi}$. Thus, $\chi\colon \mathscr{Z}_\phi\rightarrow\mathscr{Z}$ is highly proximal.

(2): Obvious by (1), for an h.p. extension is proximal and a proximal extension based on a minimal semiflow must have a unique minimal subset.

(3): Let $\mathscr{X}$ be minimal. If $\phi$ is h.p., then $Z_\phi^\natural=\{\{x\}\,|\,x\in X\}$ by (2) and obviously $\mathscr{X}\cong\mathscr{Z}_\phi^\natural$. Conversely, if $\mathscr{X}\cong\mathscr{Z}_\phi^\natural$, then $\phi=\chi$ is h.p. by (1).

(4): For necessity, assume $\phi$ is fiber-onto and open. If $t_\alpha\in S\to p\in\beta S$ and $t_\alpha z_0\to z$ in $Z$, then
$t_\alpha\phi^{-1}(z_0)=\phi^{-1}(t_\alpha z_0)\to\phi^{-1}(z)=p\diamond\phi^{-1}(z_0)$. Thus, $Z_\phi=\{\phi^{-1}(z)\colon z\in Z\}$ and $\chi$ is an isomorphism. 

Conversely, suppose $\chi\colon Z_\phi^\natural\rightarrow Z$ is 1-1. Then by (2), we have that $\chi^{-1}(z)=\{\phi^{-1}(z)\}$ for all $z\in Z$, where $\phi^{-1}(z)$ is thought of as a point in $2^{X\!,\phi}$.
Indeed, let $z\in Z$, $x\in\phi^{-1}(z)$ and $\chi^{-1}(z)=\{z^\natural\}$. Since $\mathscr{Z}$ is minimal and $\mathscr{X}$ has a dense set of a.p. points, there exists a net $t_i\in S\to p\in\beta S$ and a net of a.p. points $x_i\in\phi^{-1}(z)$ such that $t_ix_i\to x$ and $pz=z$. On the other hand, by $x_i\in z^\natural$, it follows that $x\in pz^\natural=z^\natural$. Thus, $z^\natural=\phi^{-1}(z)$.
Since $\chi$ is a homeomorphism, hence $\chi^{-1}(z)$ is continuous w.r.t. $z\in Z$ so that $\phi$ is open. To prove that $\mathscr{X}$ is $\phi$-fiber-onto, let $t\in S$ and $z\in Z$. Let $z^\natural=\phi^{-1}(z)\in Z_\phi^\natural$ and $z_1^\natural=\phi^{-1}(tz)\in Z_\phi^\natural$. Then by the $S$-equivariant of $\chi$, we have $z_1^\natural=tz^\natural$ so that $\phi^{-1}(tz)=t\diamond\phi^{-1}(z)=t\phi^{-1}(z)$. Thus, $\mathscr{X}$ is $\phi$-fiber-onto.

Finally, if $X$ is a compact metric space, then so is $2^{X,\phi}$. So $Z_\phi^\natural$, as a subspace of $2^{X,\phi}$, is a metric space as well. The proof is complete
\end{proof}
%%%%%%%%%%%%%%%%%%%%%%%%%%%%%%%%%%%%%%%%%%%%%%%%%%%%%%%%%%%%%%%%%%
\section[Open lifting via a pair of h.p. extensions]{Open lifting via a pair of h.p. extensions of semiflows}\label{SO}
This section will be devoted to proving Lemma~\ref{1.3.3} and Theorem~\ref{1.3.4} stated in $\S$\ref{s1.3}, which are contained in Theorems~\ref{3.2} and \ref{3.4}, respectively. Moreover, as a byproduct of proving Theorems~\ref{1.3.3} and \ref{1.3.4} we will solve an open question---\cite[Q.~IV.6.4b]{Wo} on some CD of semi-open extensions (Q.~\ref{3.7} $\&$ Cor.~\ref{3.9}).

If $\mathscr{X}$ and $\mathscr{Y}$ be two $S$-semiflows, then as usual $\mathscr{X}\times\mathscr{Y}$ stands for the product semiflow $S\curvearrowright_\pi X\times Y$, where $t(x,y)=(tx,ty)$ for all $t\in S$ and $(x,y)\in X\times Y$. 

As in \cite{S76, AG77, Wo, A88, DeV} in minimal flows, we can define the so-called `Auslander-Glasner commutative diagram'\,(AG-CD) for an extension $\phi\colon\mathscr{X}\rightarrow\mathscr{Z}$ of semiflows with only $\mathscr{Z}$ minimal, but with $\mathscr{X}$ having a dense set of a.p. points instead of `minimal' as follows. 

\begin{thm}[{AG-lifting; cf.~\cite{S76, AG77, Wo} or \cite[Thm.~VI.3.8]{DeV} for $\mathscr{X}$ a minimal flow}]\label{3.1}
Suppose $\phi\colon\mathscr{X}\rightarrow\mathscr{Z}$ is an extension of semiflows, where $\mathscr{Z}$ is minimal and $\mathscr{X}$ has a dense set of a.p. points. Let $\chi\colon \mathscr{Z}_\phi^\natural\rightarrow \mathscr{Z}$
be the h.p. quasifactor representation of $\mathscr{Z}$ in $\mathscr{X}$ and let
$$
X_\phi^\natural=X\vee Z_\phi^\natural=\{(x,K)\,|\,x\in K\in Z_\phi^\natural\textrm{ s.t. }\phi(x)=\chi(K)\},\quad \varrho\colon X_\phi^\natural\xrightarrow{(x,K)\mapsto x}X,\quad
\phi_\natural\colon X_\phi^\natural\xrightarrow{(x,K)\mapsto K} Z_\phi^\natural.
$$
Then there exists a `canonically determined' CD of extensions of semiflows:
$$
\leqno{\mathrm{AG}(\phi):}\qquad
\begin{CD}
\mathscr{X}@<{\varrho}<<\mathscr{X}_\phi^\natural\\
@V{\phi}VV @VV{\phi_\natural}V\\
\mathscr{Z}@<{\chi}<<\mathscr{Z}_\phi^\natural
\end{CD}
\qquad \textrm{s.t. } 
\begin{cases}
(1)\ \varrho\textrm{ and }\chi\textrm{ are h.p.};\\
(2)\ \textrm{if }\chi\textrm{ is 1-1, then so is }\varrho;\\
(3)\ \mathscr{X}_\phi^\natural\textrm{ is }\phi_\natural\textrm{-fiber-onto and }\phi_\natural\textrm{ is open}.
\end{cases}
$$
Moreover, if $X$ is a metric space and $\mathscr{X}$ is $\phi$-fiber-onto, then $\varrho$ and $\chi$ are almost 1-1 (see \cite{V70} for $\mathscr{X}$ a metric minimal flow).
\end{thm}

\begin{proof}
Clearly, $X_\phi^\natural$ is $S$-invariant closed in $X\times Z_\phi^\natural$. Since $X$ contains a dense set of a.p. points and $\mathscr{Z}$ is minimal, $\bigcup_{A\in Z_\phi^\natural}A$ is dense in $X$. Thus, $\varrho[X_\phi^\natural]=X$ and $\varrho$ is an extension of semiflows. This shows that $\mathrm{AG}(\phi)$ is a well-defined CD of extensions of semiflows, where $\mathscr{Z}$ and $\mathscr{Z}_\phi^\natural$ are minimal.

(1): Given $x\in X$ and letting $z=\phi (x)$, we have that
$$
\varrho^{-1}(x)=\{(x,K)\,|\,K\in Z_\phi^\natural\textrm{ s.t. }x\in K\}\subseteq\{x\}\times\chi^{-1}(z).
$$
Thus, $\varrho$ and $\chi$ are h.p. by Lemma~\ref{P2.3.4}. 

(2): Clearly, $\varrho$ is 1-1 restricted to each $\phi_\natural$-fiber. If $\chi$ is 1-1, then for $z_1^\natural\not=z_2^\natural$ in $Z_\phi^\natural$, then we have that $\chi(z_1^\natural)\not=\chi(z_2^\natural)$ and $\varrho[{\phi_\natural}^{-1}(z_1^\natural)]\cap\varrho[{\phi_\natural}^{-1}(z_2^\natural)]=\emptyset$ so that $\varrho$ is 1-1.

(3): By ${\phi_\natural}^{-1}(z^\natural)=K\times\{z^\natural\}$ for all $z^\natural=K\in Z_\phi^\natural$, it follows easily that ${\phi_\natural}_\textsl{ad}\colon z^\natural\mapsto{\phi_\natural}^{-1}(z^\natural)$ is continuous; and so, $\phi_\natural$ is open by Lemma~\ref{2.1.4}. 

Now for all $t\in S$ and $z^\natural=K\in Z_\phi^\natural$, since $tz^\natural=tK$ and ${\phi_\natural}^{-1}(tz^\natural)=tK\times\{tz^\natural\}=t(K\times
\{z^\natural\})$, hence $t{\phi_\natural}^{-1}(z^\natural)={\phi_\natural}^{-1}(tz^\natural)$. Thus, $\mathscr{X}_\phi^\natural$ is $\phi_\natural$-fiber-onto.

Finally, let $X$ be a metric space. Then $2^{X,\phi}$ is a metric space. Thus, $Z_\phi^\natural$ and $X_\phi^\natural$ are obviously metric spaces. 
And if, in addition, $\mathscr{X}$ is $\phi$-fiber-onto, $\varrho$ and $\chi$ are obviously almost 1-1 extensions by Lemma~\ref{P2.2.8} and Lemma~\ref{P2.3.4}. Indeed, by Lemma~\ref{P2.2.8}, we can choose a point $z_0\in Z$ such that $\phi_\textsl{ad}\colon Z\rightarrow2^X$ and $\chi_\textsl{ad}\colon Z\rightarrow2^{Z_\phi^\natural}$ both are continuous at $z_0$. Let $z\in Z$. Since $\mathscr{Z}_\phi^\natural$ is minimal (Lem.~\ref{P2.3.4}), there exists a net $t_i\in S$ with $t_iz\to z_0$ in $Z$ such that 
$t_i\phi^{-1}(z)=\phi(t_iz)\to\phi^{-1}(z_0)\in Z_\phi^\natural$ in $2^X$ and $\phi^{-1}(z_0)$ is a inclusion minimal element in $Z_\phi$. This implies that $\chi^{-1}(z_0)=\{\phi^{-1}(z_0)\}$. Thus, $\chi$ is almost 1-1. For every $x_0\in\phi^{-1}(z_0)$, since $\varrho^{-1}(x_0)=\{(x_0,\phi^{-1}(z_0))\}$, hence $\varrho$ is also almost 1-1.
The proof is complete
\end{proof}

\begin{thm}[Lem.~\ref{1.3.3}]\label{3.2}
Suppose $\phi\colon\mathscr{X}\rightarrow\mathscr{Z}$ is a fiber-onto extension of semiflows such that $\mathscr{Z}$ is minimal and $\mathscr{X}$ has a dense set of a.p. points. If $X$ is a metric space, then $\phi$ and $\phi_*$ are semi-open.
\end{thm}

\begin{proof}
    In $\mathrm{AG}(\phi)$ of Theorem~\ref{3.1}, $\phi_\natural$ is open and $\chi$ is semi-open. Thus, $\phi$ is semi-open. Further, $\phi_*$ is semi-open by Theorem~\ref{1.2.3}C. The proof is complete.
\end{proof}

However, if the factor $\mathscr{Z}$ is semi-open (Def.~\ref{1.3.1}c) instead of ``$\mathscr{X}$ is metrizable and $\phi$-fiber-onto'', then there is another generalization of Lemma~\ref{1.3.2} as follows:

\begin{thm}[{cf.~\cite[Thm.~2.5]{W85} for $S$ a group}]\label{3.3}
Let $\phi\colon\mathscr{X}\rightarrow \mathscr{Z}$ be an extension of semiflows such that $\mathscr{Z}$ is T.T. and semi-open. Then there exists a T.T. subsemiflow $S\curvearrowright X_0$ of $\mathscr{X}$ such that $\phi|_{X_0}\colon\mathscr{X}_0\rightarrow\mathscr{Z}$ is a semi-open extension.
\end{thm}

\begin{proof}
Let $\mathscr{C}$ be the collection of all closed $S$-invariant subsets of $X$ that are mapped onto $Z$ by $\phi$. Then by Zorn's Lemma there is an inclusion minimal element, say $X_0$, in $\mathscr{C}$. We shall prove that $\mathscr{X}_0$ is T.T and $\psi=\phi|_{X_0}\colon X_0\rightarrow Z$ is semi-open. Indeed, by Zorn's Lemma again, there exists a closed subset $F$ of $X_0$ that is $\psi$-irreducible (i.e., $\psi[F]=Z$ and no closed set $H\varsubsetneq F$ with $\psi[H]=Z$). Then $\psi|_F\colon F\rightarrow Z$ is semi-open. As $\mathscr{Z}$ is semi-open, it follows that for all $t\in S$, $\psi|_{tF}\colon tF\rightarrow Z$ is semi-open. Since $\psi[\overline{SF}]=Z$ and $\overline{SF}\subseteq X_0$, hence $\overline{SF}=X_0$. Moreover, $\psi\colon \overline{SF}\rightarrow Z$ is semi-open. So $\phi|_{X_0}\colon X_0\rightarrow Z$ is semi-open onto. Let $U\in\mathscr{O}(X_0)$. As $\psi[\overline{SU}]=\overline{S\psi[U]}=Z$, it follows that $\overline{SU}=X_0$. Thus, $\mathscr{X}_0$ is T.T. and this proves Theorem~\ref{3.3}.
\end{proof}

\begin{thm}[Thm.~\ref{1.3.4}]\label{3.4}
Suppose $\phi\colon\mathscr{X}\rightarrow\mathscr{Z}$ is an extension of flows such that $\mathscr{Z}$ is minimal and $\mathscr{X}$ has a dense set of a.p. points. Then $\phi$ and $\phi_*$ are semi-open.
\end{thm}

\begin{proof}
First, $\phi$ is semi-open by Lemma~\ref{1.3.2}. Next, the $\mathrm{AG}(\phi)$ implies the following CD of extensions of affine flows:
$$
\begin{CD}
\mathcal{M}^1(\mathscr{X})@<{\varrho_*}<<\mathcal{M}^1(\mathscr{X}_\phi^\natural)\\
@V{\phi_*}VV @VV{{\phi_\natural}_*}V\\
\mathcal{M}^1(\mathscr{Z})@<{\chi_*}<<\mathcal{M}^1(\mathscr{Z}_\phi^\natural)
\end{CD}.
$$

Since $\phi_\natural$ is open, hence ${\phi_\natural}_*$ is open by Theorem~\ref{2.1.5}-(1). Moreover, since $\chi$ is h.p., it follows from Lemma~\ref{P2.3.3} and Theorem~\ref{1.2.2}B that
$\chi_*$ is irreducible; and so, $\chi_*$ is semi-open. Thus, $\phi_*$ is semi-open. The proof is complete.
\end{proof}

\begin{se}\label{3.5}
    Let $\phi\colon\mathscr{X}\rightarrow\mathscr{Z}$ and $\psi\colon\mathscr{Y}\rightarrow\mathscr{Z}$ be extensions of $S$-semiflows; then the fibered product of $\phi$ and $\psi$ is defined as follows:
\begin{enumerate}
\item[] $\mathrm{R}_{\phi\psi}=\{(x,y)\in X\times Y\,|\,\phi (x)=\psi (y)\}$.
\end{enumerate}
Clearly, $\mathscr{R}_{\phi\psi}$ is a subsemiflow of $\mathscr{X}\times\mathscr{Y}$ and an extension of $\mathscr{Z}$. $\mathrm{R}_{\phi\phi}$ is usually written as $\mathrm{R}_\phi$.
\end{se}

\begin{lem}[{cf.~\cite{W85, Wo} for $\mathscr{X}$ and $\mathscr{Y}$ in the class of minimal flows}]\label{3.6}
Let there exist a CD of extensions of semiflows $\mathscr{X}$, $\mathscr{Y}$ and $\mathscr{Z}$ as follows:
$$
\begin{CD}
\mathscr{X}&@<{\quad\rho_X^{}\quad}<<&\mathscr{R}_{\phi\psi}\\
@V{\phi}VV&&@VV{\rho_Y^{}}V\\
\mathscr{Z}&@<\quad\psi\quad<<&\mathscr{Y}
\end{CD}\qquad \textrm{where} \begin{cases}\mathrm{R}_{\phi\psi}\xrightarrow{\rho_X^{}\colon (x,y)\mapsto x}X,\\ \mathrm{R}_{\phi\psi}\xrightarrow{\rho_Y^{}\colon(x,y)\mapsto y}Y.\end{cases}
$$
Then:
\begin{enumerate}[(1)]
\item If $\psi$ and $\rho_Y^{}$ are semi-open, then $\rho_X^{}$ is semi-open.
\item If $\phi$ is open, then $\rho_Y^{}$ is open.
\item If $\phi$ is open and $\psi$ is semi-open, then for every $W\in\mathscr{O}(\mathrm{R}_{\phi\psi})$ there are $U\in\mathscr{O}(X)$ and $V\in\mathscr{O}(Y)$ such that $\emptyset\not=(U\times V)\cap\mathrm{R}_{\phi\psi}\subseteq W$ and $\phi[U]=\psi[V]$.
\end{enumerate}
\end{lem}

\begin{proof}
(1): Let $W=(U\times U^\prime)\cap \mathrm{R}_{\phi\psi}$ be a basic open set in $\mathrm{R}_{\phi\psi}$. Since $\rho_Y^{}$ is semi-open, we may assume that $\rho_Y^{}[W]=U^\prime$ by shrinking $U^\prime$ if necessary. Since $\psi$ is semi-open, $V:=\textrm{int}\,\psi[U^\prime]\subseteq Z$ is nonempty. Write $U_1^\prime=\psi^{-1}[V]\cap U^\prime$ and $W_1=(U\times U_1^\prime)\cap \mathrm{R}_{\phi\psi}$ that is non-void open in $\mathrm{R}_{\phi\psi}$. Let $(x,y)\in W_1$. Let $x_n\to x$ in $X$. Then $\phi (x_n)\to\phi (x)=\psi(y)\in V$ so that there are $y_n\in U_1^\prime$ such that $(x_n,y_n)\in W_1$. Thus, $\rho_X^{}$ is semi-open.

(2): Let $(x,y)\in\mathrm{R}_{\phi\psi}$ and let $(U\times V)\cap\mathrm{R}_{\phi\psi}$ be a basic neighborhood of $(x,y)$ in $\mathrm{R}_{\phi\psi}$.
As $\phi$ is open, we may assume $\phi[U]\in\mathfrak{N}_{\phi (x)}(Z)$. Since $\psi$ is continuous and $\psi(y)=\phi (x)$, we can take $V^\prime\in\mathfrak{N}_y(Y)$ with $V^\prime\subseteq V$ such that $\psi[V^\prime]\subseteq\phi[U]$. Then $V^\prime\subseteq\rho_Y^{}[(U\times V)\cap\mathrm{R}_{\phi\psi}]$. Thus, $\rho_Y^{}$ is open.

(3): It is straightforward and we omit the details. The proof is complete.
\end{proof}

\begin{se}[{cf.~\cite[Q.~IV.6.4b]{Wo}}]\label{3.7}
Let $\phi\colon\mathscr{X}\rightarrow\mathscr{Z}$ be an extension of minimal flows. Let $\sigma_X\colon\mathscr{X}^\textrm{uhp}\rightarrow\mathscr{X}$ and 
$\sigma_Z\colon\mathscr{Z}^\textrm{uhp}\rightarrow\mathscr{Z}$ be the universal h.p. extensions, and $\phi_\textrm{uhp}\colon \mathscr{X}^\textrm{uhp}\rightarrow \mathscr{Z}^\textrm{uhp}$ 
the universal h.p. lifting of $\phi$ (cf., e.g., \cite[IV.3.10]{Wo}). Then $\phi_\textrm{uhp}$ is open such that $\phi\circ\sigma_X=\sigma_Z\circ\phi_\textrm{uhp}$ (cf.~\cite[IV]{Wo}). An open question is to ``characterize $\phi$ with $\sigma_X\times\sigma_X[\mathrm{R}_{\phi_\textrm{uhp}}]=\mathrm{R}_\phi$''. This question will be answered by Corollary~\ref{3.9}.
\end{se}

\begin{thm}\label{3.8}
We consider the CD of extensions of semiflows, where $\mathscr{X}\xleftarrow{\pi_X}\mathscr{R}_{\phi\psi}\xrightarrow{\pi_Y}\mathscr{Y}$ and $\mathscr{X}^\prime\xleftarrow{\pi_{X^\prime}}\mathscr{R}_{\phi^\prime\psi^\prime}\xrightarrow{\pi_{Y^\prime}}\mathscr{Y}^\prime$ are coordinate projections:
$$
\begin{tikzcd}
	\mathscr{X} \arrow[dd, "\phi"description] &&& \mathscr{X^{\prime}} \arrow[lll, "\sigma_1"description] \arrow[dd,dashed, "\phi^{\prime}"description,near start] &\\
	& \mathscr{R}_{\phi\psi} \arrow[ul, "\pi_X"description] \arrow[dd, "\pi_Y"description,near start] 
	&&& \mathscr{R}_{\phi^\prime \psi^{\prime}} \arrow[lll, "\sigma_1\times\sigma_2"description] \arrow[ul, "\pi_{X^\prime}"description] \arrow[dd, "\pi_{Y^\prime}"description,near start]\\
	\mathscr{Z}  &&& \mathscr{Z^{\prime}} \arrow[lll,dashed, "\tau"description] &\\
	& \mathscr{Y} \arrow[ul, "\psi"description] &&& \mathscr{Y^\prime} \arrow[lll, "\sigma_2"description] \arrow[ul,dashed, "\psi^\prime"description] &
\end{tikzcd}
\quad\quad \textrm{s.t.}
\begin{cases}
\phi^\prime\textrm{ is open},\\
\psi^\prime\textrm{ is semi-open},\\
\tau\textrm{ is irreducible}.\end{cases}
$$
Then:
\begin{enumerate}[(a)]
\item If $\sigma_1$ is semi-open, then $(\sigma_1\times\sigma_2)\mathrm{R}_{\phi^\prime\psi^\prime}=\mathrm{R}_{\phi\psi}$ if and only if  $\pi_X\colon\mathrm{R}_{\phi\psi}\rightarrow X$ is semi-open. In particular, if $\psi$ is open, then $(\sigma_1\times\sigma_2)\mathrm{R}_{\phi^\prime\psi^\prime}=\mathrm{R}_{\phi\psi}$.
\item If $\sigma_2$ is semi-open, then $(\sigma_1\times\sigma_2)\mathrm{R}_{\phi^\prime\psi^\prime}=\mathrm{R}_{\phi\psi}$ if and only if $\pi_Y\colon\mathrm{R}_{\phi\psi}\rightarrow Y$ is semi-open. In particular, if $\phi$ is open, then $(\sigma_1\times\sigma_2)\mathrm{R}_{\phi^\prime\psi^\prime}=\mathrm{R}_{\phi\psi}$.
\end{enumerate}
\end{thm}

\begin{proof}
(a): Necessity is obvious by Lemma~\ref{3.6} and $\mathrm{R}_{\phi^\prime\psi^\prime}\xrightarrow{\sigma_1\circ\pi_{X^\prime}=\pi_X\circ(\sigma_1\times\sigma_2)}X$. Now, for sufficiency, suppose that $\pi_X\colon\mathrm{R}_{\phi\psi}\rightarrow X$ is semi-open. To prove $(\sigma_1\times\sigma_2)\mathrm{R}_{\phi^\prime\psi^\prime}=\mathrm{R}_{\phi\psi}$, suppose to the contrary that $(\sigma_1\times\sigma_2)\mathrm{R}_{\phi^\prime\psi^\prime}\varsubsetneq\mathrm{R}_{\phi\psi}$. Then there exist $U\in\mathscr{O}(X)$ and $V\in\mathscr{O}(Y)$ such that
    $\emptyset\not=W:=U\times V\cap\mathrm{R}_{\phi\psi}\subseteq\mathrm{R}_{\phi\psi}\setminus(\sigma_1\times\sigma_2)\mathrm{R}_{\phi^\prime\psi^\prime}$.
    Let $U_1=\pi_X[W]$ and $V_1=\pi_Y[W]$. Since $\pi_X$ is semi-open, $\textrm{int}\,U_1\not=\emptyset$ so that $\textrm{int}\,\sigma_1^{-1}[U_1]\not=\emptyset$. As $\phi^\prime$ is open, it follows that $\phi^\prime[\sigma_1^{-1}[U_1]]\subseteq Z^\prime$ includes an open non-void subset of $Z^\prime$. Since $\tau$ is irreducible, by Lemma~\ref{P2.3.3} there exists a point $z\in\phi[U_1]=\psi[V_1]$ such that $\tau^{-1}(z)\subset\phi^\prime[\sigma_1^{-1}[U_1]]$. Now we can take a point $y\in V_1$ with $\psi (y)=z$ and then a point $y^\prime\in\sigma_2^{-1}(y)\subseteq Y^\prime$ such that $z^\prime:=\psi^\prime (y^\prime)\in\tau^{-1}(z)$. As $z^\prime\in\phi^\prime[\sigma_1^{-1}[U_1]]$, it follows that we can select a point $x^\prime\in\sigma_1^{-1}[U_1]$ with $\phi^\prime (x^\prime)=z^\prime$ and $(x^\prime,y^\prime)\in\mathrm{R}_{\phi^\prime\psi^\prime}$. Put $x=\sigma_1(x^\prime)$. Then $x\in U_1$ such that $z=\phi (x)$. Thus, $(x,y)=\sigma_1\times\sigma_2(x^\prime,y^\prime)\in W$ such that $(x,y)\in(\sigma_1\times\sigma_2)\mathrm{R}_{\phi^\prime\psi^\prime}$. This is a contradiction to  $W\cap(\sigma_1\times\sigma_2)\mathrm{R}_{\phi^\prime\psi^\prime}=\emptyset$.

(b): This may follow by a slight modification of the above proof of (a). We omit the details here. The proof is complete.
\end{proof}

\begin{cor}[{cf.~\cite[Q.~IV.6.4b]{Wo}}]\label{3.9}
Consider the following CD of extensions of minimal flows:
$$
\begin{CD}
\mathscr{X}&@<\quad \sigma\quad<<&\mathscr{X}^\prime\\
@V{\phi}VV&&@VV{\phi^\prime}V\\
\mathscr{Z}&@<\quad\tau\quad<<&\mathscr{Z}^\prime
\end{CD}
\qquad \textrm{s.t.}\begin{cases} \tau\textrm{ is h.p.},\\\phi^\prime\textrm{ is open}.\end{cases}
$$
Then $\sigma\times\sigma[\mathrm{R}_{\phi^\prime}]=\mathrm{R}_{\phi}$ if and only if $\pi_X\colon\mathrm{R}_{\phi}\rightarrow X$ is semi-open. In particular, if $\phi$ is open or $\mathrm{R}_{\phi}$ contains a dense set of a.p. points, then $\sigma\times\sigma[\mathrm{R}_{\phi^\prime}]=\mathrm{R}_{\phi}$.
\end{cor}

\begin{proof}
By Theorem~\ref{3.8} and Lemma~\ref{P2.3.3}.
\end{proof}

%%%%%%%%%%%%%%%%%%%%%%%%%%%%%%%%%%%%%%%%%%%%%%%%%%%%
\section{Quasi-separable maps and semi-openness of induced maps}\label{s4}
Recently Dai and Xie in \cite[Thm.~10C]{DX} have given a partial answer to Question~\ref{1.2.4}a using an additional condition ``densely open''. Now we will present another natural condition (Def.~\ref{4.1}) and prove that under this condition a continuous map is semi-open if and only if its induced map is semi-open (Thm.~\ref{4.5}).

\begin{se}[Quasi-separable maps]\label{4.1}
Let $f\colon X\rightarrow Y$ be a continuous onto map between compact Hausdorff spaces. Then $f$ is called \textit{quasi-separable}
    if there exists a directed set $(\Lambda,\le)$ and an inverse system $\{f_i\colon X_i\rightarrow Y\,|\,i\in\Lambda\}$ of continuous onto maps between compact metric spaces,
where $\{X_i\,|\,i\in\Lambda\}$ is an inverse system with continuous onto link maps $X_i\xleftarrow{\pi_{i,j}}X_j$ with $f_j=f_i\circ\pi_{i,j}$ for $i<j$ in $\Lambda$, such that $X=\underleftarrow{\lim}_{i\in\Lambda}\{X_i\}$ and $f=f_i\circ p_i$ (or written as $f=\underleftarrow{\lim}_{i\in\Lambda}\{f_i\}$), where $p_i\colon X\rightarrow X_i$ is the canonical projection for all $i\in\Lambda$. 
\end{se}

\begin{lem}\label{4.2}
Let $f\colon X\rightarrow Y$ be a continuous onto map, where $X$ is a compact Hausdorff space. Then $Y$ is a metrizable space if and only if $f$ is quasi-separable.
\end{lem}

\begin{proof}
Since the continuous images of compact metric spaces are metrizable, hence sufficiency is obvious. Conversely, suppose $Y$ is metrizable. Let $\Sigma(X)$ be the collection of continuous pseudo-metrics on $X$. Let $\varrho$ be the metric on $Y$. Define a partial order on $\Sigma(X)$ as follows: for any $\rho,\rho^\prime\in\Sigma(X)$, $\rho\le\rho^\prime$ iff $\rho(x,x^\prime)\le\rho^\prime(x,x^\prime)$ for all $(x,x^\prime)\in X\times X$. If $2\le n<\infty$ and $\rho_1,\dotsc,\rho_n\in\Sigma(X)$, then $\rho:=\max\{\rho_1,\dotsc,\rho_n\}\in\Sigma(X)$ such that $\rho_i\le\rho$ for $1\le i\le n$. Thus, $(\Sigma(X),\le)$ is a directed set. Now, for every $\rho\in\Sigma(X)$, define a relation on $X$ as follows:
$$
R[\rho]=\left\{(x,x^\prime)\in X\times X\,|\,\varrho(f(x),f(x^\prime))+\rho(x,x^\prime)=0\right\}.
$$
Clearly, $R[\rho]\subseteq \textrm{R}_{f}$ is a closed equivalence relation on $X$. We set $X_\rho=X/R[\rho]$, which is a compact metrizable space. Let $\lambda_\rho\colon X\rightarrow X_\rho$ and $f_\rho\colon X_\rho\rightarrow Y$ be the canonical maps. Then $f=f_\rho\circ\lambda_\rho$ for all $\rho\in\Sigma(X)$, and $R[\rho]\supseteq R[\rho^\prime]$ so that there exists a canonical link map $X_\rho\xleftarrow{\pi_{\rho,\rho^\prime}}X_{\rho^\prime}$ if $\rho\le\rho^\prime$ in $\Sigma(X)$. Thus, $\{f_\rho\colon X_\rho\rightarrow Y\,|\,\rho\in\Sigma(X)\}$ is an inverse system of continuous onto maps.
As $X=\underleftarrow{\lim}\{X_\rho\}$, it follows that $f$ is quasi-separable.
\end{proof}

\begin{lem}\label{4.3}
If $f=\underleftarrow{\lim}\{f_i \,|\,i\in\Lambda\}$ with continuous onto link maps $X_i\xleftarrow{\pi_{i,j}}X_j$, then $f$ is open (resp. semi-open) if and only if $f_i$ is open (resp. semi-open) for all $i\in\Lambda$.
\end{lem}

\begin{proof}
Necessity is obvious by $f=f_i\circ p_i$ and each $p_i\colon X\rightarrow X_i$ is a continuous onto map. Now, for sufficiency, let $U\in\mathscr{O}(X)$. By the structure of the topology of $X$, it follows that we can find some index $i\in\Lambda$ such that $p_i[U]\in\mathscr{O}(X_i)$. Then $f[U]=f_i[p_i[U]]$ is open (resp. semi-open) in $Y$. The proof is complete.
\end{proof}

We shall show that $f$ is semi-open if and only if $f_*$ is semi-open in the quasi-separable case (Thm.~\ref{4.5}). For that, we need Theorems~\ref{1.2.3}C and \ref{1.2.2}D, Lemma~\ref{4.2}, Lemma~\ref{4.3} and another lemma (Lem.~\ref{4.4}).

If $\{X_i\,|\,i\in\Lambda\}$ is an inverse system of compact metric spaces, then $\mathcal{M}^1(\prod_{i\in\Lambda}X_i)\not\cong\prod_{i\in\Lambda}\mathcal{M}^1(X_i)$ in general. 
However, we have the following with $\underleftarrow{\lim}_{i\in\Lambda}$ instead of $\prod_{i\in\Lambda}$:

\begin{lem}\label{4.4}
If $\{X_i\,|\,i\in\Lambda\}$ is an inverse system of compact Hausdorff spaces, then we have that $\mathcal{M}^1\left(\underleftarrow{\lim}_{i\in\Lambda}\{X_i\}\right)\cong\underleftarrow{\lim}_{i\in\Lambda}\left\{\mathcal{M}^1(X_i)\right\}$.
\end{lem}

\begin{proof}
Write $X_i\xleftarrow{\pi_{i,j}}X_j$ and $\pi_{i,i}=\textsl{id}_{X_i}$, for $i<j$ in $\Lambda$, for the link maps of the inverse system $\{X_i\,|\,i\in\Lambda\}$. Let $X=\underleftarrow{\lim}\{X_i\}$. Let $p_i\colon X\rightarrow X_i$ and $p_{i*}\colon\mathcal{M}^1(X)\rightarrow\mathcal{M}^1(X_i)$ be the canonical maps. As $p_i=\pi_{i,j}\circ p_j$ for $i<j$ in $\Lambda$, it follows that $p_{i*}=\pi_{i,j*}\circ p_{j*}$ so that $\{\mathcal{M}^1(X_i)\,|\,i\in\Lambda\}$ is an inverse system.
Moreover,
$\pi\colon\mathcal{M}^1(X)\rightarrow\underleftarrow{\lim}\{\mathcal{M}^1(X_i)\}$, $\mu\mapsto(p_{i*}(\mu))_{i\in\Lambda}$,
is a continuous onto map. In fact, for all $(\mu_i)_{i\in\Lambda}\in\underleftarrow{\lim}\{\mathcal{M}^1(X_i)\}$, by the finite-intersection property of compact space we may take a probability $\mu\in\bigcap_{i\in\Lambda}{p_i}_*^{-1}(\mu_i)\subseteq\mathcal{M}^1(X)$; then $\pi(\mu)=(\mu_i)_{i\in\Lambda}$.
Now let $\lambda, \mu\in\mathcal{M}^1(X)$ such that $p_{i*}(\lambda)=:\lambda_i=\mu_i:=p_{i*}(\mu)$ for all $i\in\Lambda$. We need prove that $\lambda=\mu$. By regularity of $\lambda$ and $\mu$, it is enough to show that $\lambda(K)=\mu(K)$ for all $K\in2^X$. Let $K\in2^X$ and $\varepsilon>0$. Then we can take a set $U\in\mathscr{O}(X)$ such that $\lambda(U\setminus K)+\mu(U\setminus K)<\varepsilon$ and $K\subset U$. In addition, as $i\in\Lambda$ sufficiently big, we can choose finitely many open sets, say $V_1,\dotsc,V_\ell$ in $X_i$ such that $K\subseteq p_i^{-1}(V_1)\cup\dotsm\cup p_i^{-1}(V_\ell)\subseteq U$. So by the inclusion-exclusion formula of probability or by the equality $p_i^{-1}(V_1)\cup\dotsm\cup p_i^{-1}(V_\ell)=p_i^{-1}(V_1\cup\dotsm\cup V_\ell)$, it follows that
\begin{equation*}
\lambda(K)\le\lambda(p_i^{-1}(V_1)\cup\dotsm\cup p_i^{-1}(V_\ell))
=\mu(p_i^{-1}(V_1)\cup\dotsm\cup p_i^{-1}(V_\ell))\le\mu(K)+\varepsilon.
\end{equation*}
Thus, $\lambda(K)\le\mu(K)$; and analogously, $\mu(K)\le\lambda(K)$. Then $\pi$ is 1-1 onto and $\pi$ is a homeomorphism. The proof is complete.
\end{proof}

\begin{thm}[Thm.~\ref{1.2.5}]\label{4.5}
Let $f\colon X\rightarrow Y$ be a continuous onto map between compact Hausdorff spaces. If $Y$ is metrizable, then $f$ is semi-open if and only if  $f_*\colon\mathcal{M}^1(X)\rightarrow\mathcal{M}^1(Y)$ is semi-open.
\end{thm}

\begin{proof}
Sufficiency follows from Theorem~\ref{1.2.2}D. Now, for necessity, let $Y$ be a compact metric space. Then by Lemma~\ref{4.2}, $f\colon X\rightarrow Y$ is the inverse limit of an inverse system
$\{f_i\colon X_i\rightarrow Y\,|\,i\in\Lambda\}$
of continuous onto maps with $X_i$ compact metrizable. Let $p_i\colon X\rightarrow X_i$ be the canonical map. Then $f=f_i\circ p_i$ for all $i\in\Lambda$.
Thus, $f_i\colon X_i\rightarrow Y$ is semi-open by Lemma~\ref{4.3}. We have then concluded a CD of continuous onto maps:
$$
\begin{tikzcd}
	\mathcal{M}^1(X)
	\arrow[drr, "{p_i}_*"description]
	\arrow[drrr, "{p_j}_*"description, bend left=10]
	\arrow[dd, "f_*" description] &&&& \\
	& \cdots & \mathcal{M}^1(X_i) \arrow[l] \arrow[lld, "{f_i}_*"description] & \mathcal{M}^1(X_j) \arrow[l, "{\pi_{i,j}}_*"'] \arrow[llld, "{f_j}_*"description, bend left=10]
	&\cdots \arrow[l]& \underleftarrow{\lim}\{\mathcal{M}^1(X_i)\} \arrow[l] \\
	\mathcal{M}^1(Y)&&&&
\end{tikzcd}
$$
So by Theorem~\ref{1.2.3}C, it follows that ${f_i}_*\colon \mathcal{M}^1(X_i)\rightarrow\mathcal{M}^1(Y)$ is semi-open for all $i\in\Lambda$ so that by Lemmas~\ref{4.3} and \ref{4.4}, $f_*\colon\mathcal{M}^1(X)\rightarrow\mathcal{M}^1(Y)$ is semi-open. The proof is complete.
\end{proof}

\begin{rem}\label{4.6}
    Let $\xi\in\mathcal{M}^1(X\times Y)$, where $X$, $Y$ are compact metric spaces. Let $\mu\in\mathcal{M}^1(X)$ be the projection of $\xi$ onto $X$. Then there is a random probability measure $\nu_\cdot\colon X\rightarrow\mathcal{M}^1(Y)$ such that $\xi=\int\nu_xd\mu(x)$. If $\nu_x\equiv\nu$ for all $x\in X$, then $\xi=\mu\otimes\nu$. But $\nu_x\not\equiv\nu\ \forall x\in X$ in general. Thus, $\mathcal{M}^1(X\times Y)\not\cong\mathcal{M}^1(X)\times\mathcal{M}^1(Y)$ in general.
\end{rem}

\begin{rem}[Suggested by the referee]\label{4.7}
  Let $f=\underleftarrow{\lim}\{f_i\}\colon X\rightarrow Y$ where $\{f_i\colon X_i\rightarrow Y\,|\,i\in\Lambda\}$ is an inverse system of irreducible continuous onto mappings with continuous onto link maps $X_i\xleftarrow{\pi_{i,j}}X_j$ for $i<j$ in $\Lambda$. Then $f$ is irreducible.
\end{rem}

Consequently, if $f\colon X\rightarrow Y$ is the inverse limit of an inverse system $\{f_i\colon X_i\rightarrow Y\,|\,i\in\Lambda\}$ of almost 1-1 mappings, then $f$ is irreducible; but $f$ is generally not almost 1-1.

\begin{exa}\label{4.8}
    Let $X=Y\times\{1\}\sqcup Y\times\{2\}$ be Ellis' two-circle space where $Y=\mathbb{S}$ is the unit circle in $\mathbb{C}$ (Ex.~\ref{2.2.2}F). Let $f\colon X\rightarrow Y$ be the projection $(a,i)\mapsto a$.
    Then $f$ is 2-1 and irreducible.
    Let $(\Lambda,\le)$ be the directed set of all continuous pseudo-metrics on $X$ as in Proof of Lemma~\ref{4.2}. Given $d\in\Lambda$, let $X_d=X/R[d]$ and let $X\xrightarrow{p_d}X_d\xrightarrow{f_d}Y$ be the canonically induced maps as in Proof of Lemma~\ref{4.2}, where 
    $R[d]=\left\{(x,x^\prime)\in X\times X\colon |f(x)-f(x^\prime)|+d(x,x^\prime)=0\right\}$ 
    is a closed equivalence subrelation of $\mathrm{R}_f$ on $X$.
    Clearly, $f_d\colon X_d\rightarrow Y$ is irreducible by $f=f_d\circ p_d$; and moreover, $f=\underleftarrow{\lim}_{d\in\Lambda}\left\{f_d\right\}$ by Lemma~\ref{4.2} and $X\cong\underleftarrow{\lim}_{d\in\Lambda}\{X_d\}$ ($\because X$ is a compact Hausdorff space and $\bigcap_{d\in\Lambda}R[d]=\varDelta_X$). However, since $X_d$, for each $d\in\Lambda$, is a compact metric space, hence $f_d$ is almost 1-1 by Lemma~\ref{P2.2.8}. Thus, there exists a 2-1 irreducible map $f$ that is the inverse limit of an inverse system of almost 1-1 mappings.
\end{exa}

%%%%%%%%%%%%%%%%%%%%%%%%%%%%%%%%%%%%%%%%%%%%%%%%%%%%
\section{Quasi-separable extensions of minimal semiflows}\label{s5}
This section will be devoted to proving Theorem~\ref{1.3.5} (Thm.~\ref{5.2.2}) based on Lemma~\ref{1.3.3} and considering quasi-separable minimal weakly mixing flows (Thm.~\ref{Q5.10} $\&$ Thm.~\ref{Q5.11}). 

\subsection{Basic notions}
In what follows, let $S$ be a topological monoid, not necessarily discrete. Then, for any semiflow $\mathscr{X}=S\curvearrowright_\pi X$, we require the phase transformation $\pi\colon S\times X\rightarrow X$, $(s,x)\mapsto sx$, is a jointly continuous mapping.
Let $\mathscr{X}$ be a semiflow with phase semigroup $S$. Then:

\begin{sse}\label{5.1.1}
If every point of $X$ is a.p. (cf.~Def.~\ref{1.3.1}b), then $\mathscr{X}$ is termed a \textit{pointwise a.p. semiflow}. 
\end{sse}

\begin{sse}\label{5.1.2}
We say that $\mathscr{X}$ is \textit{algebraically transitive} (A.T.) if $Sx=X\ \forall x\in X$; If $\mathscr{X}\times\mathscr{X}$ is T.T. (cf.~Def.~\ref{1.3.1}a), then $\mathscr{X}$ is termed \textit{weakly mixing}.
\end{sse}

\begin{sse}\label{5.1.3}
Let    
\begin{enumerate}
    \item[] $P(\mathscr{X})=\left\{(x,x^\prime)\in X\times X\,|\,\overline{S(x,x^\prime)}\cap\Delta_X\not=\emptyset\right\}$,
    \end{enumerate}
    which is called the \textit{proximal relation} on $\mathscr{X}$. For all $x\in X$ let
    \begin{enumerate}
    \item[] $P[x]=\left\{x^\prime\in X\,|\,(x,x^\prime)\in P(\mathscr{X}\right\}$,
    \end{enumerate}
    which is called the \textit{proximal cell} at $x$ of $\mathscr{X}$.
    $\mathscr{X}$ is called a \textit{proximal flow} if $P(\mathscr{X})=X\times X$. If $P(\mathscr{X})=\Delta_X$, then $\mathscr{X}$ is said to be \textit{distal}.  See, e.g., \cite{E69, G76, B79, Wo, A88, DeV, CD}.
\end{sse}

\begin{sse}\label{5.1.4} 
If every minimal proximal flow with phase group $S$ is a singleton, then $S$ is termed \textit{strongly amenable}\,(\cite[$\S$II.3]{G76}). For example, the compact extension of a nilpotent group is strongly amenable\,(see, e.g., \cite[Thm.~II.3.4]{G76} or \cite[Prop.~1.4]{MW74}).
\end{sse}

\begin{sse}\label{5.1.5}
We say that $\mathscr{X}$ is a \textit{Bron\v{s}te\v{\i}n semiflow} (\textsl{B}-semiflow; cf.~\cite{V77}) if $\mathscr{X}\times\mathscr{X}$ has a dense set of a.p. points. It turns out that if $\mathscr{X}$ is a minimal flow with strongly amenable phase group, then $\mathscr{X}$ is a \textsl{B}-flow\,(see, e.g., \cite[Prop.~X.1.3]{G76}).
\end{sse}

\begin{srem}\label{5.1.6}
Let $\mathscr{X}$ be an $S$-semiflow. Clearly, $\pi_*\colon S\times\mathcal{M}^1(X)\rightarrow\mathcal{M}^1(X)$, $(t,\mu)\mapsto t\mu=t_*\mu$ is separately continuous such that $e\mu=\mu$, $(st)\mu=s(t\mu)$ for all $\mu\in\mathcal{M}^1(X)$ and $s,t\in S$. Thus, if $S$ is locally compact Hausdorff group, then $(t,\mu)\mapsto t\mu$ is jointly continuous so that $S\curvearrowright_{\pi_*}\!\mathcal{M}^1(X)$ is a flow by Ellis's joint continuity theorem (\cite{E57}).

\begin{enumerate}[$\bullet$]
\item In fact, $(t,\mu)\mapsto t\mu$ is still jointly continuous and $S\curvearrowright_{\pi_*}\!\mathcal{M}^1(X)$, denoted $\mathcal{M}^1(\mathscr{X})$, is an affine semiflow in general.
\end{enumerate}
    \begin{proof}
    First, if $t_n\to t$ in $S$, then $t_nx\to tx$ in $X$ uniformly for $x\in X$. For otherwise, there is an $\varepsilon\in\mathscr{U}_X$ the uniformity structure of $X$ and a subnet $\{t_{n^\prime}\}$ of $\{t_n\}$ and $x_{n^\prime}\in X\to x^\prime$ such that $(t_{n^\prime}x_{n^\prime}, tx_{n^\prime})\notin\varepsilon$. So $(tx^\prime,tx^\prime)\notin\varepsilon$, a contradiction.
    Second, if $f\in C(X)$ and $t_n\to t$ in $S$, then $\|ft_n-ft\|_\infty\to 0$ by the first assertion.
    Last, let $\mu_n\to\mu$ in $\mathcal{M}^1(X)$ and $t_n\to t$ in $S$. By the uniform bounded principle or the resonance theorem, it follows that $\|\mu_n\|\cdot\|ft_n-ft\|_\infty\to 0$ for all $f\in C(X)$. Thus, for all $f\in C(X)$,
    $$    
    \lim|\mu_n(ft_n)-\mu(ft)|\le\lim(\|\mu_n\|\cdot\|ft_n-ft\|_\infty+|\mu_n(ft)-\mu(ft)|)=0.
    $$
    This shows that $t_n\mu_n\to t\mu$ in $\mathcal{M}^1(X)$. Whence $(t,\mu)\mapsto t\mu$ is jointly continuous.
    \end{proof}
\end{srem}

\subsection{Quasi-separable extensions}\label{s5.2}
Recall that a quasi-separable mapping (Def.~\ref{4.1}) does not involve any dynamics and it is not necessarily to be semi-open. Next, we shall introduce a ``dynamical quasi-separable mapping'' and prove that it is automatically semi-open under the fiber-onto condition.

\begin{sse}[Quasi-separable extensions]\label{5.2.1}
Let $\phi\colon \mathscr{X}\rightarrow \mathscr{Y}$ be an extension of $S$-semiflows. $\phi$ is called \textit{quasi-separable} (cf.~\cite{E69, K71, K72, E73}), if there exists a directed set $(\Lambda,\le)$ and an inverse system of extensions of semiflows
$\{\phi_i\colon\mathscr{X}_i\rightarrow\mathscr{Y}\,|\,i\in\Lambda\}$,
where $\{\mathscr{X}_i\}$ is an inverse system of metrizable $S$-semiflows, such that $\mathscr{X}=\underleftarrow{\lim}\{\mathscr{X}_i\}$ and $\phi=\phi_i\circ\rho_i$ (or written as $\phi=\underleftarrow{\lim}\{\phi_i\}$), where $\rho_i\colon \mathscr{X}\rightarrow \mathscr{X}_i$ is the canonical projection for all $i\in\Lambda$ such that $\rho_i=\pi_{i,j}\circ\rho_j$ and $\mathscr{X}_i\xleftarrow{\pi_{i,j}}\mathscr{X}_j$ is the connection extension for all $i<j$ in $\Lambda$. In the special case $\mathscr{Y}$ is a singleton, $\mathscr{X}$ is called a \textit{quasi-separable semiflow}.
\end{sse}

It should be mentioned that for any non-metrizable $\mathscr{X}$, $\phi\colon X\rightarrow\{y\}$ is a quasi-separable mapping by Lemma~\ref{4.2}; however, $\phi\colon\mathscr{X}\rightarrow\{y\}$ is not necessarily to be a quasi-separable extension (cf.~Rem.~\ref{5.2.3}). Thus, the concept of quasi-separable semiflow is non-trivial.

\begin{sthm}[Thm.~\ref{1.3.5}]\label{5.2.2}
Let $\phi\colon\mathscr{X}\rightarrow\mathscr{Y}$ be a fiber-onto quasi-separable extension of minimal $S$-semiflows. Then $\phi$ and $\phi_*$ are semi-open.
\end{sthm}

\begin{proof}
In view of Theorem~\ref{1.2.5} (or Thm.~\ref{4.5}), it suffices to prove that $\phi$ is semi-open. For that, let
$$\phi=\underleftarrow{\lim}\left\{\mathscr{X}_i\xrightarrow{\phi_i}\mathscr{Y}\,|\,i\in\Lambda\right\}$$
as in Def.~\ref{5.2.1}, where each $\mathscr{X}_i$ is a minimal metrizable $S$-semiflow.
Since $\mathscr{X}$ is $\phi$-fiber-onto, it follows by $\phi=\phi_i\circ\rho_i$ that each $\mathscr{X}_i$ is also $\phi_i$-fiber-onto. Thus, by Lemma~\ref{1.3.3} (or Thm.~\ref{3.2}), $\phi_i$ is semi-open for each $i\in\Lambda$. Finally, $\phi$ is semi-open by Lemma~\ref{4.3}.
\end{proof}

\begin{srem}[{cf.~\cite{E68, K72}, \cite[Lem.~2.1]{K71} or \cite[Thm.~I.1.7]{Wo} for $\mathscr{X}$ to be point-transitive by using $S$-subalgebra of $C(\beta S)$ with $S$ a discrete group}]\label{5.2.3}
\textit{If $\mathscr{X}$ has $\sigma$-compact phase group, then
$\mathscr{X}$ is a quasi-separable flow.} Subsequently, a flow having separable locally compact phase group is quasi-separable.
\end{srem}

\begin{proof}
Let $S=\bigcup_{n=1}^\infty K_n$, where each $K_n$, $n\in\mathbb{N}$, is a compact subset of $S$. Let $\rho\in\Sigma(X)$ and $R[\rho]=\{(x,x^\prime)\in X\times X\,|\,\rho(tx,tx^\prime)=0\ \forall t\in S\}$. Then $R[\rho]$ is an invariant closed equivalence relation on $X$. Set $X_i=X/R[i]$ for all $i\in\Sigma(X)$, where the quotient space $X_i$ is metrizable via the metric 
$$
d_i([x]_{R[i]},[x^\prime]_{R[i]})=\sum_{n=1}^\infty\frac{1\wedge\max\{\rho(tx,tx^\prime)\colon t\in K_n\}}{2^n}\quad \forall x,x^\prime\in X.
$$
Clearly, $\{\mathscr{X}_i\,|\,i\in\Sigma(X)\}$ is an inverse system of metrizable $S$-flows. Let $\lambda_i\colon \mathscr{X}\rightarrow\mathscr{X}_i$ be the canonical map for all $i\in\Sigma(X)$. Then
$\lambda\colon\mathscr{X}\rightarrow\underleftarrow{\lim}\{\mathscr{X}_i\}$, given by $x\mapsto(\lambda_ix)_{i\in\Sigma(X)}$,
is 1-1. To prove that $\lambda$ is onto, let $(x_i)_{i\in\Sigma(X)}\in\underleftarrow{\lim}\{\mathscr{X}_i\}$ be arbitrary. Set $[x_i]=\{x\in X\,|\,\lambda_ix=x_i\}$, which is a closed nonempty subset of $X$. For any $i_1,\dotsc,i_n\in\Sigma(X)$, there is some $j\in\Sigma(X)$ with $i_1\le j, \dotsc, i_n\le j$.
Then $[x_j]\subseteq[x_{i_1}]\cap\dotsm\cap[x_{i_n}]$. Thus, $\bigcap_{i\in\Sigma(X)}[x_i]\not=\emptyset$. Take $x\in\bigcap_{i\in\Sigma(X)}[x_i]$. Then $\lambda x=(x_i)_{i\in\Sigma(X)}$. Thus, $\lambda$ is onto and $\mathscr{X}\cong\underleftarrow{\lim}\{\mathscr{X}_i\}$ is a quasi-separable flow.
\end{proof}

However, if $\phi\colon\mathscr{X}\rightarrow\mathscr{Y}$ is an extension of flows with $\sigma$-compact phase group $S$, where $Y$ need not be metrizable, we don't know whether or not $\phi$ is quasi-separable.

\begin{sq}\label{5.2.4}
Let $\phi\colon\mathscr{X}\rightarrow\mathscr{Y}$ be an extension of minimal \textit{semiflows}, where $\mathscr{X}$ is either metrizable or $\phi$-fiber-onto non-quasi-separable. {\it Is $\phi$ semi-open?}
\end{sq}

\subsection{Quasi-separable minimal flows}\label{s5.3}
Quasi-separable minimal flows are special quasi-separable extensions of minimal flows. In this subsection we will digress from ``semi-openness'' to study the structure of minimal weakly mixing quasi-separable flows.

\begin{slem}\label{Q5.7}
Let $\mathscr{X}$ be a pointwise a.p. T.T. flow. If $\mathscr{X}$ is quasi-separable, then $\mathscr{X}$ is minimal.
\end{slem}

\begin{proof}
There exists an inverse system $\{\mathscr{X}_i,\pi_{i,j}\}$ of metrizable flows such that $\mathscr{X}=\underleftarrow{\lim}\{\mathscr{X}_i\}$. Then $\mathscr{X}_i$ is pointwise a.p. T.T. with $X_i$ a compact metrizable space. Thus, $\mathscr{X}_i$ is minimal so that $\mathscr{X}$ is minimal. The proof is complete.
\end{proof}

\begin{slem}\label{Q5.8}
Let $\{\mathscr{X}_i\,|\,i\in\Lambda\}$ be an inverse system of flows with links $\{\pi_{i,j}\colon i,j\in\Lambda\}$ and let $\mathscr{X}=\underleftarrow{\lim}\{\mathscr{X}_i\}$. Then:
\begin{enumerate}[(1)]
\item $\{\mathscr{X}_i\times\mathscr{X}_i, \pi_{i,j}\times\pi_{i,j}\,|\,i\in\Lambda\}$ is an inverse system of flows such that $\mathscr{X}\times\mathscr{X}=\underleftarrow{\lim}\{\mathscr{X}_i\times\mathscr{X}_i\}$.
\item $\mathscr{X}$ is T.T. if and only if $\mathscr{X}_i$ is T.T. for all $i\in\Lambda$.
\item $\mathscr{X}$ is weak-mixing if and only if $\mathscr{X}_i$ is weak-mixing for all $i\in\Lambda$ (cf.~\cite[Lem.~2.5]{K71} for $\mathscr{X}$ minimal and $S$ abelian).
\item If $\mathscr{X}$ is minimal, then $\overline{P(\mathscr{X})}=X\times X$ if and only if $\overline{P(\mathscr{X}_i)}=X_i\times X_i$ for all $i\in\Lambda$.
\end{enumerate}
\end{slem}

\begin{proof}
(1): Obvious.

(2): Necessity is evident. Now for sufficiency, let $U$, $V\in\mathscr{O}(X)$. Then there exists some $i\in\Lambda$ and there are $U_i$, $V_i\in \mathscr{O}(X_i)$ such that $p_i^{-1}[U_i]\subseteq U$ and $p_i^{-1}[V_i]\subseteq V$, where $p_i\colon\mathscr{X}\rightarrow\mathscr{X}_i$ is the canonical map. Since $\mathscr{X}_i$ is T.T., $tU_i\cap V_i\not=\emptyset$ for some $t\in S$. So $tU\cap V\not=\emptyset$. Thus, $\mathscr{X}$ is T.T.

(3): By (1) and (2).

(4): Since $\mathscr{X}$ is a minimal flow and $p_i\colon\mathscr{X}\rightarrow\mathscr{X}_i$ is onto, hence $p_i\times p_i[P(\mathscr{X})]=P(\mathscr{X}_i)$ for all $i\in\Lambda$. Thus, necessity is obvious. Now suppose $\overline{P(\mathscr{X}_i)}=X_i\times X_i$ for all $i\in\Lambda$. Let $U$, $V\in\mathscr{O}(X)$. As $X=\underleftarrow{\lim}\{X_i\}$, it follows that there exist $i\in\Lambda$ and $U_i$, $V_i\in \mathscr{O}(X_i)$ such that $p_i^{-1}[U_i]\subseteq U$ and $p_i^{-1}[V_i]\subseteq V$. Further, there is a pair $(x_i,y_i)\in U_i\times V_i\cap P(\mathscr{X}_i)$. Clearly, there is a pair $(x,y)\in U\times V\cap P(\mathscr{X})$ with $p_i\times p_i(x,y)=(x_i,y_i)$. Thus, $P(\mathscr{X})$ is dense in $X\times X$.
The proof is complete.
\end{proof}

\begin{slem}[{cf.~\cite[Prop.~II.2.1]{G76}}]\label{Q5.9}
Let $\mathscr{X}$ be a flow such that for every $m\ge2$ and for all $A, U_1, \dotsc, U_m\in\mathscr{O}(X)$, we have that
$(U_1\times\dotsm\times U_m)\cap S[A\times\dotsm\times A]\not=\emptyset$.
If $\mathscr{Y}$ is a minimal flow, then $\mathscr{X}\times\mathscr{Y}$ is a T.T. flow.
\end{slem}

\begin{proof}
Let $U, A\in\mathscr{O}(X)$ and $V, B\in \mathscr{O}(Y)$. We need prove that $(U\times V)\cap S^{-1}[A\times B]\not=\emptyset$.
Since $\mathscr{Y}$ is a minimal flow, there are finitely many elements $t_1, \dotsc,t_m\in S$ with $t_1^{-1}[V]\cup\dotsm\cup t_m^{-1}[V]=Y$. Then there are points $a_1, \dotsc,a_m\in A$ and $s\in S$ such that
$s(a_1, \dotsc,a_m)\in t_1^{-1}[U]\times\dotsm\times t_m^{-1}[U]$. Take $b\in B$. As $sb\in Y$, it follows that $sb\in t_k^{-1}[V]$ for some $1\le k\le m$. Thus,
$s(a_k, b)\in t_k^{-1}[U\times V]$ and $t_ks(a_k,b)\in U\times V$. The proof is complete.
\end{proof}

\begin{sthm}[{Wu's problem~\cite[Prob.~7, p.~518]{AG68}: \textit{$(1)\Rightarrow(3)$?}}]\label{Q5.10}
Let $\mathscr{X}$ be a minimal quasi-separable flow satisfying one of the following conditions:
\item $\mathrm{C}_1$. $\mathscr{X}$ is a \textsl{B}-flow;
\item $\mathrm{C}_2$. $\mathscr{X}$ admits a regular Borel probability measure.
\item Then the following are pairwise equivalent:
\begin{enumerate}[(1)]
\item $\mathscr{X}$ has no non-trivial distal factor.
\item $P[x]$ is dense in $X$ for all $x\in X$.
\item $P(\mathscr{X})$ is dense in $X\times X$.
\item $\mathscr{X}$ is weakly mixing.
\end{enumerate}
(\textbf{Note.} See, e.g.,~\cite[Thm.~9.13]{A88} and \cite{DeV} for the case that $X$ is a compact metric space.)
\end{sthm}

\begin{proof}
Let
$R\!P(\mathscr{X})=\{(x,x^\prime)\in X\times X\,|\,\overline{S[U\times V]}\cap\Delta_X\not=\emptyset\ \forall\, U\times V\in\mathfrak{N}_{(x,x^\prime)}(X\times X)\}$; and let
$$
U(\mathscr{X})=\{(x,x^\prime)\in X\times X\,|\,\exists x_i^\prime\to x^\prime, t_i\in S\textrm{ s.t. }t_i(x,x_i^\prime)\to(x,x)\}.
$$
Then, under condition $\mathrm{C}_1$\,(cf.~\cite[Thm.~2.7.6]{V77}) or $\mathrm{C}_2$\,(cf.~\cite{Mc78}), we have that $R\!P(\mathscr{X})=U(\mathscr{X})$ is an invariant closed equivalence relation. Thus,
$$
R\!P[x]={\bigcap}_{\varepsilon\in\mathscr{U}_X}\overline{S^{-1}[\varepsilon[x]]}\quad \textrm{and}\quad P[x]={\bigcap}_{\varepsilon\in\mathscr{U}_X}S^{-1}[\varepsilon[x]]
$$
for all $x\in X$. Let $\mathscr{X}=\underleftarrow{\lim}\{\mathscr{X}_i\}$ where $\mathscr{X}_i$, $i\in\Lambda$, are minimal metrizable flows. Let $\rho_i\colon\mathscr{X}\rightarrow\mathscr{X}_i$ be the canonical maps.

$(1)\Rightarrow(2)$: First, by Furstenberg's structure theorem, it follows that $R\!P(\mathscr{X})=X\times X$. As $\rho_i\times\rho_i[R\!P(\mathscr{X})]=R\!P(\mathscr{X}_i)$, it follows that $R\!P(\mathscr{X}_i)=X_i\times X_i=U(\mathscr{X}_i)$ for all $i\in\Lambda$. Hence $R\!P[x_i]=U[x_i]=X_i$ for all $x_i\in X_i$ and $i\in\Lambda$. Since $X_i$ is metrizable, $\mathscr{U}_{X_i}$ has a countable basis. Thus, $\overline{P[x_i]}=X_i$ for all $x_i\in X_i$ and all $i\in\Lambda$. Now let $x\in X$ and set $x_i=\rho_i(x)$. Let $U\in \mathscr{O}(X)$. Then there exist $i\in\Lambda$ and $U_i\in\mathscr{O}(X_i)$ with $\rho_i^{-1}[U_i]\subseteq U$. We can take a point $y_i\in U_i$ such that $y_i\in P[x_i]$. Further, there exists a point $y\in X$ such that $y\in P[x]$ and $\rho_i(y)=y_i$. So $y\in\rho_i^{-1}(y_i)\subseteq U$. This shows that $\overline{P[x]}=X$.

$(2)\Rightarrow(3)\Rightarrow(1)$: Obvious.

\item Consequently, $(1)\Leftrightarrow(2)\Leftrightarrow(3)$. It remains to prove that $(1)\Leftrightarrow(4)$.

$(4)\Rightarrow(1)$: Obvious.

$(2)\Rightarrow(4)$: Let $m\ge2$. (2) implies that $(U_1\times\dotsm\times U_m)\cap S[A\times\dotsm\times A]\not=\emptyset$ for all $A$, $U_1, \dotsc, U_m\in \mathscr{O}(X)$. Thus, $\mathscr{X}\times\mathscr{X}$ is T.T. so $\mathscr{X}$ weakly mixing by Lemma~\ref{Q5.9}. The proof is complete.
\end{proof}

Note that $\mathrm{C}_2$ need not imply $\mathrm{C}_1$ in Theorem~\ref{Q5.10}. Indeed, Furstenberg's example (\cite[II.5.5]{G76}) says that there exists a solvable group $S$ such that there is a non-trivial proximal (so non-B) $S$-flow. However, we have no example of ``$\mathrm{C}_1\not\Rightarrow\mathrm{C}_2$'' at hands.

\begin{sthm}\label{Q5.11}
Let $\mathscr{X}$ be a pointwise a.p., weakly mixing, quasi-separable, non-trivial flow with strongly amenable phase group. Then:
\begin{enumerate}[(1)]
\item $P(\mathscr{X})$ is not an equivalence relation\,(cf.~\cite[Prop.~2.3]{K71} for $\mathscr{X}$ minimal $\&$ $S$ abelian).
\item There exists no closed proximal cell for $\mathscr{X}$. (Consequently, $\mathscr{X}$ is not point distal.)
\end{enumerate}
\end{sthm}

\begin{proof}
First $\mathscr{X}$ is minimal by Lemma~\ref{Q5.7}. Since $\mathscr{X}$ is quasi-separable, there exists an inverse system $\{\mathscr{X}_i\,|\,i\in\Lambda\}$ of minimal metrizable flows such that $\mathscr{X}=\underleftarrow{\lim}\{\mathscr{X}_i\}$. Let $p_i\colon\mathscr{X}\rightarrow\mathscr{X}_i$ be the canonical map for all $i\in\Lambda$. Then $\mathscr{X}_i$ is a weak-mixing minimal metrizable \textsl{B}-flow for all $i\in\Lambda$.

(1): To prove that $P(\mathscr{X})$ is not an equivalence relation, suppose to the contrary that $P(\mathscr{X})$ is an equivalence relation. Then $P(\mathscr{X}_i)=p_i\times p_i[P(\mathscr{X})]$ is also an equivalence relation. Thus, by \cite[Cor.~2.14]{MW74}, it follows that $X_i=P\circ P[x_i]=P[x_i]$ for all $x_i\in X_i$ and all $i\in\Lambda$. (In fact, for $x_i,x_i^\prime\in X_i$, $P[x_i]$ and $P[x_i^\prime]$ are residual in $X_i$; then $P[x_i]\cap P[x_i^\prime]\not=\emptyset$ implies that $(x_i,x_i^\prime)\in P(\mathscr{X}_i)$, so $P[x_i]=X_i$.) Whence $\mathscr{X}_i$ is proximal so that $X_i$ is a singleton for all $i\in\Lambda$, for $S$ is strongly amenable. This shows that $\mathscr{X}$ is a trivial flow, contrary to that $\mathscr{X}$ is non-trivial.

(2): Suppose to the contrary that $P[x]$ is a closed set for some point $x\in X$. Then it is easy to verify that $P[x_i]$ is closed in $X_i$ for all $i\in\Lambda$, where $x_i=p_i(x)$. By Theorem~\ref{Q5.10}, it follows that $P[x_i]=X_i$ for all $i\in\Lambda$. Thus, $\mathscr{X}_i$ is proximal and $X_i$ is a singleton. Then $\mathscr{X}$ is a singleton, a contradiction.
The proof is complete.
\end{proof}

\begin{srem}\label{Q5.12}
The ``strongly amenable'' condition is crucial in Theorem~\ref{Q5.11}. For example, let $S=\mathrm{SL}(2,\mathbb{R})$ the topological group of real $2\times 2$ matrices with determinant $1$. Let $X=\mathbb{P}^1$ be the projective line, i.e., the set of lines through the origin of the plane. $S$ acts naturally on $\mathbb{P}^1$ (sending lines into lines), and this action is A.T. Then $\mathscr{X}$ is minimal proximal (then not a B-flow)\,(see \cite[II.5.6]{G76}), and $\mathscr{X}$ is weakly mixing~(see \cite[Cor.~II.2.2]{G76} or Lem.~\ref{Q5.9}).
Moreover, $\mathscr{X}$ admits no invariant Borel probability measures. Otherwise, suppose $\mu$ be an invariant Borel probability measure; as $\mathrm{SL}(2,\mathbb{R})$ includes the rotations, it follows that $\mu$ is the Lebesque measure; this contradicts the proximality. Thus, neither $\mathrm{C}_1$ nor $\mathrm{C}_2$ is a necessary condition for Theorem~\ref{Q5.10}.
\end{srem}

%%%%%%%%%%%%%%%%%%%%%%%%%%%%%%%%%%%%%%%%%%%%%%%%%%%%%%%%%%%%%%%%%%
\section{Pontryagin's open-mapping theorem}\label{s6}
In this section we shall consider the semi-openness of some canonical mappings following from the theory of topological groups (Thm.~\ref{6.1} and Cor.~\ref{6.3}) and an application in topological groups (Thm.~\ref{6.4}).

\begin{thm}[Thm.~\ref{1.4.3}]\label{6.1}
Let $G$ and $H$ be left-topological groups, such that $G$ is Lindel\"{o}f  quasi-regular and $H$ is Baire. If $f\colon G\rightarrow H$ is a locally closed continuous onto $G$-equivariant function, then $f$ is a semi-open mapping.
\end{thm}

\begin{proof}
Let $S=G$ be equipped with the discrete topology. First $S\curvearrowright G$ is a minimal $S$-flow. Let $S\curvearrowright_\pi H$ be defined by $\pi_t\colon H\rightarrow H$, $h\mapsto \tau(t)h$, for all $t\in S$, where $\tau\colon G\rightarrow H$ is given as in Def.~\ref{1.4.2}. Then $S\curvearrowright_\pi H$ is also a minimal $S$-flow and $f\colon\mathscr{G}\rightarrow\mathscr{H}$ is a locally closed extension of $S$-flows (not necessarily compact Hausdorff). So $f$ is semi-open by Theorem~\ref{P2.2.1}.
The proof is complete.
\end{proof}

Notice here that if $G$ is a topological group and $f$ is a homomorphism in Theorem~\ref{6.1}, then for $U\in\mathfrak{N}_e(G)$ we can take $V\in\mathfrak{N}_e(G)$ with $V=V^{-1}$ and $V^2\subseteq U$. This implies that $f[U]\in\mathfrak{N}_e(H)$; and so, $f$ is open (Pontryagin's open-mapping theorem---Thm.~\ref{1.4.1}).

\begin{cor}[Thm.~\ref{1.4.4}]\label{6.2}
Let $G$ and $H$ be left-topological groups, such that $G$ is locally compact Lindel\"{o}f pseudo-metrizable and $H$ is Hausdorff Baire. If $f\colon G\rightarrow H$ is a continuous onto $G$-equivariant function, then $f$ is an open mapping.
\end{cor}

\begin{proof}
    By Theorem~\ref{6.1}, $f$ is a semi-open continuous onto mapping. Then by Lemma~\ref{P2.2.7}, it follows that there exists a point $x_0\in G$ such that $f[U]\in\mathfrak{N}_{f(x_0)}(H)$ for all $U\in\mathfrak{N}_{x_0}(G)$. Finally, by the equivariance of $f$, we have for every $x\in G$ that $f[U]\in\mathfrak{N}_{f(x)}(H)$ for all $U\in\mathfrak{N}_{x}(G)$. Thus, $f$ is an open mapping.
\end{proof}

\begin{cor}\label{6.3}
Let $G$ be a right-topological group on a locally compact Lindel\"{o}f Hausdorff space and $g\in G$. If $L_g\colon G\rightarrow G$ is continuous, then $L_g$ is semi-open; and moreover, $L_g$ is open if $G$ is developable.
\end{cor}

\begin{proof}
    By Theorem~\ref{6.1} and Lemma~\ref{P2.2.7}.
\end{proof}

\begin{thm}\label{6.4}
Let $G$ be a right-topological group on a locally compact, Lindel\"{o}f, developable, Hausdorff space. If $g\in G$ such that
$L_g\colon G\rightarrow G$ is continuous, then $L_{g^{-1}}\colon G\rightarrow G$ is also continuous.
\end{thm}

\begin{proof}
By Corollary~\ref{6.3}, it follows that $L_g\colon G\rightarrow G$ is a homeomorphism; and so, $L_{g^{-1}}\colon G\rightarrow G$ is also continuous.
\end{proof}

%%%%%%%%%%%%%%%%%%%%%%%%%%%%%%%%%%%%%%%%%%%%%%%%%%%%%%%%%%%%%%%%%%%
\section*{Acknowledgements}
\noindent This research was partially supported by the National Natural Science Foundation of China (Grant No. 12271245). The authors would like to thank the referee for her/his constructive comments.


\end{document}